\documentclass[11pt]{article}
\usepackage{amssymb}
\newtheorem{thm}{Claim}[subsection]
\newtheorem{propose}[thm]{Proposition}
\newtheorem{lemma}[thm]{Lemma}
\newtheorem{cor}[thm]{Scholium}
\newtheorem{question}[thm]{Problem}
\newtheorem{defn}[thm]{Definition}
\newtheorem{conj}[thm]{Conjecture}
\newtheorem{remark}[thm]{Remark}

\renewcommand{\d}{\mbox{\LARGE $\cdot $}}
\newcommand{\Xs}{X_{\d}}            

\newcommand{\Hom}{{\rm Hom}\,}      
\newcommand{\Ext}{{\rm Ext}\,}      

\renewcommand{\P}{\mbox{$\mathbb P$}}   
\newcommand{\Q}{\mbox{$\mathbb Q$}}     
\newcommand{\A}{\mbox{$\mathbb A$}}     
\newcommand{\C}{\mbox{$\mathbb C$}}     
\newcommand{\Z}{\mbox{$\mathbb Z$}}     
\newcommand{\HH}{\mbox{$\mathbb H$}}    
\newcommand{\R}{\mbox{$\mathbb R$}}     
\newcommand{\bPic}{\mbox{$\bf Pic$}}

\newcommand{\im}{{\rm im}\,}        
\newcommand{\gr}{{\rm gr}\,}        
\newcommand{\Pic}{{\rm Pic}\,}     
\newcommand{\NS}  {{\rm NS}\,}      

\newcommand{\rank}{{\rm rank}\,}    

\newcommand{\by}[1]{\stackrel{#1}{\rightarrow}}
\newcommand{\longby}[1]{\stackrel{#1}{\longrightarrow}}

\newcommand{\df}{\mbox{\,$\stackrel{\rm def}{=}$}\,}
\newcommand{\ie}{{\it i.e.\/},\ }
\newcommand{\cf}{{\it cf.\/}\ }
\newcommand{\eg}{{\it e.g.\/},\ }

\renewcommand{\bar}{\overline}
\newcommand{\into}{\hookrightarrow}
\newcommand{\implies}{\mbox{$\Rightarrow$}}
\newcommand{\veq}{\mbox{\large $\parallel$}}  
\newcommand{\sZ}{\mbox{\scriptsize{$\Z$}}}   
\newcommand{\sC}{\mbox{\scriptsize{$\C$}}}   
\newcommand{\sQ}{\mbox{\scriptsize{$\Q$}}}   

\newcommand{\limdir}[1]{\mathop{\rm
lim}_{\buildrel\longrightarrow\over{#1}}}

\newcommand{\onto}{\mbox{$\to\!\!\!\!\to$}}

\newcommand{\boxtensor}{\def\boxtimesten{\Box\kern-7.59pt\raise1.2pt
\hbox{$\times$} }}                                  

\newcounter{elno}                   

\newenvironment{proof}{{\bf Proof}:\quad
                     }{{\hfill$\odot$\\}}

\newcommand{\cA}{{\cal A}}

\newcommand{\cD}{{\cal D}}

\newcommand{\cF}{{\cal F}}

\newcommand{\cH}{{\cal H}}

\newcommand{\cK}{{\cal K}}

\newcommand{\cO}{{\cal O}}

\newcommand{\cW}{{\cal W}}

\newcommand{\cZ}{{\cal Z}}
\catcode`\@=11
\def\opn#1#2{\def#1{\mathop{\kern0pt\fam0#2}\nolimits}} 
\def\underrightarrow{\mathpalette\underrightarrow@}
\def\underrightarrow@#1#2{\vtop{\ialign{$##$\cr
 \hfil#1#2\hfil\cr\noalign{\nointerlineskip}%
 #1{-}\mkern-6mu\cleaders\hbox{$#1\mkern-2mu{-}\mkern-2mu$}\hfill
 \mkern-6mu{\to}\cr}}}

\def\underleftarrow{\mathpalette\underleftarrow@}
\def\underleftarrow@#1#2{\vtop{\ialign{$##$\cr
\hfil#1#2\hfil\cr\noalign{\nointerlineskip}#1{\leftarrow}\mkern-6mu
 \cleaders\hbox{$#1\mkern-2mu{-}\mkern-2mu$}\hfill
 \mkern-6mu{-}\cr}}}
\def\:{\colon}
\let\oldtilde=\tilde

\def\indextil#1{\lower2pt\hbox{$\textstyle{\oldtilde{\raise2pt%
\hbox{$\scriptstyle{#1}$}}}$}}
\def\pnt{{\raise1.1pt\hbox{$\textstyle.$}}}
\let\amp@rs@nd@\relax
\newdimen\ex@
\newdimen\bigaw@
\newdimen\minaw@
\newdimen\minCDaw@
\newif\ifCD@
\def\minCDarrowwidth#1{\minCDaw@#1}

\def\@CD{\def\A##1A##2A{\llap{$\vcenter{\hbox
 {$\scriptstyle##1$}}$}\Big\uparrow\rlap{$\vcenter{\hbox{%
$\scriptstyle##2$}}$}&&}%
\def\V##1V##2V{\llap{$\vcenter{\hbox
 {$\scriptstyle##1$}}$}\Big\downarrow\rlap{$\vcenter{\hbox{%
$\scriptstyle##2$}}$}&&}%
\def\={&\hskip.5em\mathrel
 {\vbox{\hrule width\minCDaw@\vskip3\ex@\hrule width
 \minCDaw@}}\hskip.5em&}%
\def\verteq{\Big\Vert&&}%
\def\noarr{&&}%
\def\vspace##1{\noalign{\vskip##1\relax}}\relax\let\amp@rs@nd@&\iffalse}\fi
 \CD@true\vcenter\bgroup\relax\let\\=\cr\iffalse}\fi\tabskip\z@skip\baselineskip20\ex@
 &\hfill$\m@th##$\hfill\cr}
\def\@endCD{\cr\egroup\egroup}
\def\>#1>#2>{\amp@rs@nd@\setbox\z@\hbox{$\scriptstyle
 \;{#1}\;\;$}\setbox\@ne\hbox{$\scriptstyle\;{#2}\;\;$}\setbox\tw@
 \hbox{$#2$}\ifCD@
 \global\bigaw@\minCDaw@\else\global\bigaw@\minaw@\fi
 \ifdim\wd\z@>\bigaw@\global\bigaw@\wd\z@\fi
 \ifdim\wd\@ne>\bigaw@\global\bigaw@\wd\@ne\fi
 \ifCD@\hskip.5em\fi
 \ifdim\wd\tw@>\z@
 \mathrel{\mathop{\hbox to\bigaw@{\rightarrowfill}}\limits^{#1}_{#2}}\else
 \mathrel{\mathop{\hbox to\bigaw@{\rightarrowfill}}\limits^{#1}}\fi
 \ifCD@\hskip.5em\fi\amp@rs@nd@}
\def\<#1<#2<{\amp@rs@nd@\setbox\z@\hbox{$\scriptstyle
 \;\;{#1}\;$}\setbox\@ne\hbox{$\scriptstyle\;\;{#2}\;$}\setbox\tw@
 \hbox{$#2$}\ifCD@
 \global\bigaw@\minCDaw@\else\global\bigaw@\minaw@\fi
 \ifdim\wd\z@>\bigaw@\global\bigaw@\wd\z@\fi
 \ifdim\wd\@ne>\bigaw@\global\bigaw@\wd\@ne\fi
 \ifCD@\hskip.5em\fi
 \ifdim\wd\tw@>\z@
 \mathrel{\mathop{\hbox to\bigaw@{\leftarrowfill}}\limits^{#1}_{#2}}\else
 \mathrel{\mathop{\hbox to\bigaw@{\leftarrowfill}}\limits^{#1}}\fi
 \ifCD@\hskip.5em\fi\amp@rs@nd@}
\def\@CDS{\def\A##1A##2A{\llap{$\vcenter{\hbox
 {$\scriptstyle##1$}}$}\Big\uparrow\rlap{$\vcenter{\hbox{%
$\scriptstyle##2$}}$}&}%
\def\V##1V##2V{\llap{$\vcenter{\hbox
{$\scriptstyle##1$}}$}\Big\downarrow\rlap{$\vcenter{\hbox{%
$\scriptstyle##2$}}$}&}%
\def\={&\hskip.5em\mathrel
 {\vbox{\hrule width\minCDaw@\vskip3\ex@\hrule width
 \minCDaw@}}\hskip.5em&}
\def\verteq{\Big\Vert&}
\def\novarr{&}
\def\noharr{&&}
\def\SE##1E##2E{\slantedarrow(0,18)(4,-3){##1}{##2}&}
\def\SW##1W##2W{\slantedarrow(24,18)(-4,-3){##1}{##2}&}
\def\NE##1E##2E{\slantedarrow(0,0)(4,3){##1}{##2}&}
\def\NW##1W##2W{\slantedarrow(24,0)(-4,3){##1}{##2}&}
\def\slantedarrow(##1)(##2)##3##4{%
\thinlines\unitlength1pt\lower 6.5pt\hbox{\begin{picture}(24,18)%
\put(##1){\vector(##2){24}}%
\put(0,8){$\scriptstyle##3$}%
\put(20,8){$\scriptstyle##4$}%
\end{picture}}}
\def\vspace##1{\noalign{\vskip##1\relax}}\relax\let\amp@rs@nd@&\iffalse}\fi
 \CD@true\vcenter\bgroup\relax\let\\=\cr\iffalse}\fi\tabskip\z@skip\baselineskip20\ex@
 &\hfill$\m@th##$\hfill\cr}
\def\@endCDS{\cr\egroup\egroup}
\newdimen\TriCDarrw@
\newif\ifTriV@

\def\@TriCDV{\TriV@true\def\TriCDpos@{6}\@TriCD}
\def\@TriCDA{\TriV@false\def\TriCDpos@{10}\@TriCD}
\def\@TriCD#1#2#3#4#5#6{%
\setbox0\hbox{$\ifTriV@#6\else#1\fi$}
\TriCDarrw@=\wd0 \advance\TriCDarrw@ 24pt
\advance\TriCDarrw@ -1em
\def\SE##1E##2E{\slantedarrow(0,18)(2,-3){##1}{##2}&}
\def\SW##1W##2W{\slantedarrow(12,18)(-2,-3){##1}{##2}&}
\def\NE##1E##2E{\slantedarrow(0,0)(2,3){##1}{##2}&}
\def\NW##1W##2W{\slantedarrow(12,0)(-2,3){##1}{##2}&}
\def\slantedarrow(##1)(##2)##3##4{\thinlines\unitlength1pt
\lower 6.5pt\hbox{\begin{picture}(12,18)%
\put(##1){\vector(##2){12}}%
\put(-4,\TriCDpos@){$\scriptstyle##3$}%
\put(12,\TriCDpos@){$\scriptstyle##4$}%
\end{picture}}}
\def\={\mathrel {\vbox{\hrule
   width\TriCDarrw@\vskip3\ex@\hrule width
   \TriCDarrw@}}}
\def\>##1>>{\setbox\z@\hbox{$\scriptstyle
 \;{##1}\;\;$}\global\bigaw@\TriCDarrw@
 \ifdim\wd\z@>\bigaw@\global\bigaw@\wd\z@\fi
 \hskip.5em
 \mathrel{\mathop{\hbox to \TriCDarrw@ {\rightarrowfill}}\limits^{##1}}
 \hskip.5em}
\def\<##1<<{\setbox\z@\hbox{$\scriptstyle
 \;{##1}\;\;$}\global\bigaw@\TriCDarrw@
 \ifdim\wd\z@>\bigaw@\global\bigaw@\wd\z@\fi
 \mathrel{\mathop{\hbox to\bigaw@{\leftarrowfill}}\limits^{##1}}
 }
 \CD@true\vcenter\bgroup\relax\let\\=\cr\iffalse}\fi
 \tabskip\z@skip\baselineskip20\ex@
 \ifTriV@
 \halign\bgroup
 &\hfill$\m@th##$\hfill\cr
\else
 \halign\bgroup
 &\hfill$\m@th##$\hfill\cr
\fi}
\def\@endTriCD{\egroup}

\title{\sc On algebraic 1-motives related to Hodge cycles}

\author{by Luca {\sc Barbieri-Viale}}

\date{}

\begin{document}

\pagestyle{headings}

\maketitle

\begin{abstract}
The goal of this paper is to introduce {\em Hodge 1-motives} of
algebraic varieties and to state a corresponding cohomological Grothendieck-Hodge
conjecture, generalizing the classical Hodge conjecture to arbitrarily
singular proper schemes.

We also construct generalized cycle class maps from the $\cK$-cohomology
groups $H^{p+i}(\cK_p)$ to the sub-quotiens $W_{2p}H^{2p+i}/W_{2p-2}$
given by the weight filtration. However, in general, the image of this cycle
map (as well as the image of the canonical map from motivic cohomology) is
strictly smaller than the rational part of the Hodge  filtration $F^p$ on
$H^{2p+i}$.
\end{abstract}

\tableofcontents

\section{Introduction}
Let $X$ be an algebraic $\C$-scheme.  The singular cohomology groups
$H^*(X,\Z (\cdot))$ carry a mixed Hodge structure, see \cite[III]{D}.
Deligne theory of 1-motives (see \cite[III]{D}) is an algebraic framework
in order to deal with {\em some} mixed Hodge structures extracted from
$H^*(X,\Z (\cdot))$,
\ie those having non-zero Hodge numbers in the set $\{(0,0), (0,-1), (-1,0),
(-1,-1)\}$. Therefore, these cohomological invariants of algebraic
varieties would be algebraically defined as 1-motives over arbitrary base
fields or schemes.  Note that a general theory of mixed motives can be
regarded as an algebraic framework in order to deal with {\em all} mixed
Hodge structures $H^*(X,\Z (\cdot))$.

A 1-motive $M$ over a scheme $S$ is given by an $S$-homomorphism of group
schemes $$M = [L\by{u} G]$$ where $G$ is an extension of an abelian scheme
$A$ by a torus $T$ over $S$, and the group scheme $L$ is, locally for the
\'etale topology on $S$, isomorphic to a given finitely-generated free
abelian group. There are Hodge, De Rham and $\ell$-adic realizations (see
\cite[III]{D} and \cite{DM}).

If $X$ is a smooth proper $\C$-scheme then $H^i(X,\Z(j))$ is pure of weight
$i-2j$.  If $i=2p$ is even a natural 1-motive would be given by the lattice
of Hodge cycles in $H^{2p}(X,\Z(p))$, \ie of those integral cohomology
classes (modulo torsion) which are of type $(0,0)$.  Classical Hodge
conjecture claims that (over $\Q$) such a 1-motive would be obtained from classes of
algebraic cycles on $X$ only.  For $i=2p+1$ odd the
1-motive corresponding to $H^{2p+1}(X,\Z(p+1))$ is given by the abelian
variety associated to the largest sub-Hodge structure whose types are
$(-1,0)$ or $(0,-1)$.  Grothendieck-Hodge conjecture characterize (over $\Q$)
this sub-Hodge structure as the coniveau $\geq p$ sub-space, \ie the abelian
variety as the algebraic part of the intermediate jacobian.

Grothendieck-Hodge conjectures are concerned with the quest of an algebraic
definition for the named 1-motives.  In fact, the usual Hodge conjecture
can be reformulated by saying that the Hodge realization of the
algebraically defined $\Q$-vector space of codimension $p$ algebraic cycles
modulo numerical (or homological) equivalence is the 1-motivic part of
$H^{2p}(X,\Q(p))$.  Moreover, the 1-motivic part of $H^{2p+1}(X,\Q(p+1))$
would be the Hodge realization of the isogeny class of the universal
regular quotient.

The main task of this paper is to define {\em Hodge 1-motives}\, of singular
varieties and to state a corresponding cohomological Grothendieck-Hodge
conjecture, by dealing with their Hodge realizations.

\subsection*{A short survey of the subject}
The classical Hodge conjecture along with a tantalizing overview can be
found in \cite{DH}. Recall that Grothendieck corrected the general Hodge
conjecture in \cite{GH}. The book of Lewis \cite{LW} is a very good compendium of
methods and results.

Recall that Jannsen \cite{JA} formulated an homological version of the
Hodge conjecture for singular varieties.  Moreover, Bloch in a letter to
Jannsen (see the Appendix~A in \cite{JA} \cf Section 5), gave a counterexample
to a naive cohomological Hodge conjecture for curves on a singular 3-fold.
However, in the same letter, Bloch was guessing that the Hodge conjecture for
divisors, \ie $F^1\cap H^2(X,\Z)$ is generated by $c_1$ of line bundles on
$X$, holds true in the singular setting ``because one has the exponential''.
Anyways, jointly with V. Srinivas, we gave a counterexample to this claim
and questioned a reformulation of the Hodge conjecture for divisors in
\cite{BS} by restricting to Zariski locally trivial cohomology classes,
\ie let $L^pH^*(X, \Z)$ be the filtration induced by the Leray spectral
sequence along the canonical continuous map $X_{\rm an}\to X_{\rm Zar}$,
is $F^1\cap L^1H^2(X, \Z)$ given by $c_1$ of line bundles on $X$ ?
Still, this reformulation doesn't hold in general, \eg see \cite{BiS} where
it is also proved for $X$ normal.

 From the work of Carlson (see \cite{CA} and \cite{C}) and the theory of
Albanese and Picard 1-motives \cite{BSAP} it now appears that the theory of
1-motives is a natural setting for a formulation of a cohomological version
of the Hodge conjectures for singular varieties.  For example, $F^1\cap
H^2(X,\Z)$ is simply given by $c_1$ of simplicial line bundles on a smooth
proper hypercovering $\pi : X_{\d}\to X$ via universal cohomological
descent $\pi^* :H^2(X,\Z)\cong H^2(X_{\d},\Z)$.  This N\'eron-Severi group
$\NS (\Xs) \cong F^1\cap H^2(X,\Z)$ is actually independent of the choice
of the smooth simplicial scheme.  Furthermore, $\NS (\Xs)$ admits an
algebraic definition as the quotient of the simplicial Picard group scheme
$\bPic_{\Xs}$ by its connected component of the identity (\cf \cite{BSAP}).
However, the 1-motivic part of $H^2(X,\Z)$ is still larger than $\NS
(\Xs)$. Therefore the largest algebraic part of $H^2(X,\Z)$ will be
detected from a honest 1-motive only (see \cite{BRS}).

Note that the rank of the usual $\NS (X)$ (= the image of $\Pic (X)$ in
$H^2(X,\Z)$) is actually smaller than $\NS (\Xs)$, in general.  Moreover
$F^1\cap W_0H^2(X,\Q) =0$ thus $\NS (\Xs)_{\sQ}$ is naturally a subspace of
$H^2(X,\Q)/W_0$.

Considering the Leray filtration $L^pH^{2p}(X, \Q)$ a natural
question formulated in \cite{BiS} is if $F^p\cap L^pH^{2p}(X, \Q)$ will be
given by higher Chern classes. However, this would not be true without
some extra hypothesis on $X$ and does not tell enough about the algebraic
part of all $H^{2p}(X, \Q)$.

\subsection*{An outline of the conjectural picture}
Let $X$ be a proper integral $\C$-scheme.  Let $H \df H^{2p+i}(X)$ be our
mixed Hodge structure on $H^{2p+i}(X, \Z)/ {\rm (torsion)}$ for a fixed
pair of integers $p\geq 0$ and $i\in\Z$.  First remark that we always have
an extension $$0\to\gr^W_{2p-1}H \to W_{2p}H/W_{2p-2}H\to \gr^W_{2p}H\to 0.
$$ An extension always defines an extension class map $$e^p: H^{p,p}_{\sZ}
\to J^p(H)$$ which is not, in general, a 1-motive.  In fact, $J^p(H)$ is a
complex torus which is not an abelian variety, in general.  Recall that
Carlson \cite{CA} studied abstract extensions of Hodge structures showing
their geometric content for low-dimensional varieties.  For higher
dimensional schemes consider the largest abelian subvariety $A^p(H)$ of the
torus $J^p(H)$.  Denote $H^p(H)$ the group of Hodge cycles, that is the
preimage in $H^{p,p}_{\sZ}$ of $A^p(H)$ under the extension class map.
Note that an example due to Srinivas shows that the group $H^p(H)$ of Hodge cycles
in this sense can be strictly smaller than $H^{p,p}_{\sZ}$ (see
Section~5.2).

Define the {\it Hodge 1-motive}\, of the mixed Hodge structure $H$ the so
obtained 1-motive $$e^p: H^{p}(H)\to A^p(H).$$
Conversely, Deligne's theory of 1-motives \cite{D} grant us of a mixed
Hodge structure $H^h$ corresponding to this 1-motive (see Section~2.2 for
details). According to Deligne's philosophy of 1-motives there should be
an algebraically defined 1-motive whose Hodge realization is $H^h$.
The algebraic definition (see Section~2.1) is predictable {\it via}\,
Grothendieck-Hodge conjectures and Bloch-Beilinson motivic world as
follows.

Assume $X$ smooth.  Consider the filtration $F^*_a$ on the Chow group
$CH^p(X)$ given by $F^0_a = CH^p(X)$, $F^1_a = CH^p(X)_{\rm alg}$ the
sub-group of cycles algebraically equivalent to zero and $F^2_a=$ the
kernel of the Abel-Jacobi map.  Thus the graded pieces are $\gr^0_{F_a} =
NS^p(X)$, the N\'eron-Severi group, and $\gr^1_{F_a} = J^p_a(X)=$ the group
of $\C$-points of an abelian subvariety of the intermediate jacobian.  We
then get an extension $$0 \to J^p_a(X) \to CH^p(X)/F^2_a \to NS^p(X)\to
0.$$ Note that Grothendieck-Hodge conjecture claims that
$J^p_a(X)=A^p(H^{2p-1}(X))$ up to isogeny, \ie $J^p_a(X)$ is the largest
abelian subvariety of the intermediate jacobian.

If $X$ is not smooth then let $\Xs$ be a smooth proper simplicial scheme
along with $\pi : \Xs \to X$, a universal cohomological descent morphism
(\cf \cite{GRO}).  In zero characteristic, such $\Xs$ was firstly provided
by the construction of hypercoverings in
\cite{D}, then by that of cubical hyperresolutions in \cite{GN} where the
dimensions of the components are bounded or by the method
of hyperenvelopes given in \cite{GS}. In \cite{DJ} such a simplicial
scheme is provided in positive characteristics.

The above extension, given by the filtration $F^*_a$ on the Chow groups of
each component of the so obtained simplicial scheme $\Xs$, yields a short
exact sequence of complexes.  Let $(NS^p)^{\bullet}$ and
$(J^p_a)^{\bullet}$ denote such complexes.  By taking homology groups we
then get boundary maps $$\lambda^i_a : H^i((NS^p)^{\bullet}) \to
H^{i+1}((J^p_a)^{\bullet}).$$ We conjecture that the boundary map
$\lambda^i_a$ behave well with respect to the extension class map $e^p$
yielding a motivic cycle class map, \ie the following diagram
$$\begin{array}{ccc}
H^{i}((NS^p)^{\bullet})&\by{\lambda^i_a}& H^{i+1}((J^p_a)^{\bullet})\\
\downarrow & &\downarrow\\
H^{2p+i}(X)^{p,p} & \by{e^p} & J^p(H^{2p+i}(X))
\end{array}$$
commutes. Note that all maps in the square are canonically defined.
We guess that the image 1-motive (up to isogeny!) is the above Hodge
1-motive of the mixed
Hodge structure (see Conjecture \ref{MHC} for a full statement).
In fact, we may expect $J^p_a$ would be obtained as the universal regular
quotient of $CH^p(X)_{\rm alg}$ and that the filtration $F^*_a$ would be
induced by the motivic filtration conjectured by Bloch, Murre and
Beilinson. Accordingly we can sketch an algebraic definition of such Hodge
1-motives (see Section~2.1).

If $X$ is singular one is then puzzled by the role of $$\Hom_{\rm MHS}
(\Z(-p), H^{2p+i}(X))$$ where the $\Hom$ is taken in the abelian category
of mixed Hodge structures.  That is the integral part of $F^p$ (the Hodge
filtration on $H^{2p+i}(X,\C)$) which is contained in the kernel of the
extension class map $e^p$ above. Note that $F^p\cap W_{2p-2}H^{2p+i}(X,\Q)=
0$ here. In the smooth case, such a target is
usually reached by algebraic cycles.  In order to obtain cycle class maps
we may use local higher Chern classes and edge maps in coniveau spectral
sequences (see \cite{BO} and \cite{BV2}).

In the singular case, we show that such edge maps can be recovered by
weight arguments. In order to do this we define {\em Zariski sheaves}\, of
mixed Hodge structures, obtaining {\em infinite dimensional}\, mixed Hodge
structures on their cohomology (see Section~3). The main example is given
by the Zariski sheaf $\cH^*_X$ associated to the presheaf $U\subset X
\mapsto H^*(U)$ of mixed Hodge structures.  On the smooth simplicial scheme
$\Xs$ we also have a simplicial sheaf $\cH_{\Xs}^*$ of mixed Hodge
structures.  Since $\pi : \Xs \to X$ yields $H^*(X)\cong \HH^{*}(X_{\d})$
we then obtain a local-to-global spectral sequence $$L^{p,q}_2 =
\HH^p(X_{\d}, \cH_{\Xs}^q)\implies H^{p+q}(X)$$ in the category of infinite
dimensional mixed Hodge structures.  The sheaf $\cH_{\Xs}^q$ has weights
$\leq 2q$ and the same holds for its cohomology.  There is an edge map (see
Section~4) $$s\ell^{p+i}: \HH^{p+i}(X_{\d}, \cH_{\Xs}^p)/W_{2p-2} \to
W_{2p}H^{2p+i}(X)/W_{2p-2}.$$ We expect that the image of $s\ell^{p+i}$ is
the mixed Hodge structure $H^{2p+i}(X)^h$ corresponding to the Hodge
1-motive.

Moreover, consider $\cK$-cohomology groups $\HH^*(X_{\d},
\cK_p)$ where $\cK_p$ are the simplicial sheaves associated to Quillen's
higher $K$-theory.  Recall that there are local higher Chern classes
$$c_p : \cK_p \to \cH_{\Xs}^p(p)$$ for each $p\geq 0$.
We thus obtain a generalized cycle class map
$$c\ell^{p+i} : \HH^{p+i}(X_{\d}, \cK_p)_{\sQ} \to
W_{2p}H^{2p+i}(X, \Q)/W_{2p-2}.$$ However the image of $c\ell^{p+i}$ is not
$F^p\cap H^{2p+i}(X, \Q)$, \ie the rational part of the Hodge filtration can be
larger (see Section~5.1 where Bloch's counterexample is explained). The same
applies to the canonical map $H^{2p+i}_m(X, \Q (p)) \to H^{2p+i}(X, \Q (p))$ from
motivic cohomology.

\subsection*{Towards Hodge mixed motives}
Any reasonable theory of mixed motives would include the theory of
1-motives, \ie it would be a fully faithful functor from the $\Q$-linear
category of 1-motives to that of mixed motives.  This is the case of the
triangulated category of geometrical motives introduced by Voevodsky (see
\cite[3.4]{V}, \cf \cite{LM} and \cite{LI}).  Hanamura's construction (see
\cite{HA} and \cite{HAM}) doesn't apparently provide such a property as yet.

As remarked by Grothendieck \cite[\S 2]{GH} and Deligne \cite[\S 5]{DH} the
Hodge conjecture yields nice properties of the Hodge realization of pure
motives, \ie the usual Hodge conjecture means that the Hodge realization
functor is fully faithful. It would be interesting to investigate such a
property in the mixed case, \eg if this formulation of the Hodge
conjecture provide such a property of mixed motives.

We remark that M.~Saito recently observed (see \cite[2.5 (ii)]{MS} and
\cite{MSC}) that the canonical functor from arithmetic Hodge structures to
mixed Hodge structures is not full. Even if the Hodge realization factors
through arithmetic  Hodge structures, this non-fullness  doesn't imply the
non-fullness of the  Hodge realization of mixed motives (as noticed by
M.~Saito as well).

However, the first natural attempt to go further with Hodge mixed motives
is to provide an intrinsic definition of such objects internally.  In fact,
since 1-motives provide mixed motives we may claim that such Hodge mixed
motives exist and would be naturally defined over any field or base scheme.

\subsection*{Acknowledgements}
Papers happen because many people co-operate with one another. Often the
author is just one of a whole group of people who pitch in, and that was
the case here.

Gratitude is due therefore, firstly, to V.~Srinivas: he has been
exceedingly generous to me with his time, shared enthusiasms over the
years and told quite a few tricks.

A huge intellectual debt is due to P.~Deligne who also interceded with
advice that helped on several portions of the manuscript.

I am grateful also to S.~Bloch, O.~Gabber, H.~Gillet, M.~Hanamura, U.~Jannsen,
M.~Levine, J.~D.~Lewis, J.~Murre, A.~Rosenschon, M.~Saito, C.~Soul\'e and
V.~Voevodsky for discussions on some matters treated herein.

This research was carried out with smooth efficiency thanks to several
foundations. I like to mention Tata Institute of Fundamental Research and
Institut Henri Poincar\'e for their support and hospitality.

\section{Filtrations on Chow groups}
Following the general framework of mixed motives (\eg see \cite{DM},
\cite{LM} and \cite{JM} for a full overview) we may expect the following
picture for non-singular algebraic varieties over a field $k$
(algebraically closed of characteristic zero for simplicity).

Let $X$ be a smooth proper $k$-scheme. Bloch, Beilinson and Murre (\cf
\cite{JM}, \cite{JF} and \cite{MU3}) conjectured the existence of a
finite filtration $F_m^*$ on Chow groups $CH^p(X)\df \cZ^p(X)/\equiv_{\rm rat}$
of codimension $p$ cycles modulo rational equivalence such that
\begin{itemize}
\item $F_m^0CH^p(X)=CH^p(X)$,
\item $F_m^1CH^p(X)$ is given by $CH^p(X)_{\rm hom}$, \ie by the sub-group
of those codimension $p$ cycles which are homologically equivalent to zero
for some Weil cohomology theory,
\item $F_m^*CH^p(X)$ should be functorial and compatible with the
intersection pairing, and
\item this filtration should be motivic, \eg $\gr_{F_m}^iCH^p(X)_{\sQ}$
depends only on the Grothendieck motive $h^{2p-i}(X).$
\end{itemize}

\subsection{Regular homomorphisms}
Consider the sub-group $CH^p(X)_{\rm alg}$ of those cycles in $CH^p(X)$
which are algebraically equivalent to zero, \ie $CH^p(X)_{\rm alg}\df \ker
(CH^p(X)\to NS^p(X)).$ Denote $CH^p(X)_{\rm ab}$ the sub-group of
$CH^p(X)_{\rm alg}$ of those cycles which are abelian equivalent to zero,
\ie  $CH^p(X)_{\rm ab}$ is the intersection of all kernels of regular
homomorphisms from $CH^p(X)_{\rm alg}$ to abelian varieties.

Assume the existence of a universal regular homomorphism $\rho^p
:CH^p(X)_{\rm alg}\to A^p_{X/k}(k)$ to (the group of $k$-points) of an
abelian variety $A^p_{X/k}$ defined over the base field $k$ (\cf \cite{LB}).
This is actually proved for $p=1, 2, \dim(X)$ (see \cite{MU1}).

Thus $CH^p(X)_{\rm alg}\subseteq CH^p(X)_{\rm hom}$
and there would be an induced functorial filtration $F_a^*$ on
$CH^p(X)$ such that
\begin{itemize}
\item $F_a^0CH^p(X) = CH^p(X)$,
\item $F_a^1CH^p(X) = CH^p(X)_{\rm alg}$ is the sub-group of cycles
algebraically
equivalent to zero, and
\item $F_a^2CH^p(X) = CH^p(X)_{\rm ab}$, \ie is the kernel of the universal
regular homomorphism $\rho^p$ defined above.
\end{itemize}
Note that the existence of the abelian variety $A^p_{X/k}$ is not
explicitly mentioned in the context of mixed motives but is a rather
natural property after the case $k=\C$.

For $X$ smooth and proper over $\C$ one obtains that the motivic filtration
is such that {\it (i)}\, $F_m^1CH^p(X)=CH^p(X)_{\rm hom}$ is the sub-group
of cycles whose cycle class in $H^{2p}(X,\Z(p))$ is zero, and {\it (ii)}\,
$F_m^2CH^p(X)$ is contained in the kernel of the Abel-Jacobi map
$CH^p(X)_{\rm hom}\to J^p(X)$ and $F_m^2CH^p(X)\cap CH^p(X)_{\rm alg}$ is the
kernel of the Abel-Jacobi map $CH^p(X)_{\rm alg}\to J^p(X)$.
In this case, $CH^p(X)_{\rm ab}$ will be the kernel of the Abel-Jacobi map
$CH^p(X)_{\rm alg}\to J^p(X)$, \ie {\it (iii)}\, $A^p_{X/\sC}$ is the
algebraic part of the intermediate jacobian.  Is well known that the image
$J^p_a(X)$ of $CH^p(X)_{\rm alg}$ into $J^p(X)$ yields a sub-torus of $J^p(X)$
which is an abelian variety: moreover, is known to be universal for $p=1, 2,
\dim(X)$ (see \cite{MU1} and \cite{MU2}).

In the following, for the sake of simplicity, the reader can indeed assume
that $k =\C$ and {\it (i)--(iii)}\, are satisfied by the first two steps of
the filtration $F^i_m$.  In fact, for $k=\C$, S. Saito has obtained (up to
torsion!)  such a result (see \cite[Prop.  2.1]{SA}, \cf \cite{LF} and
\cite{JF}).  Moreover, in the following, the reader could also avoid
reference to the motivic filtration by dealing with the first
two steps of the ``algebraic'' filtration $F^i_a$ defined above.

\subsection{Extensions}
Let $X$ be smooth and proper over $k$. For our purposes just consider the
following extension
\begin{equation}\label{cyclext}
0\to \gr^1_{F_m}CH^p(X) \to CH^p(X)/F_m^2 \to \gr^0_{F_m}CH^p(X) \to 0
\end{equation}
Note that $A^p_{X/k}(k)$ is contained in $\gr^1_{F_m}CH^p(X)$ (since
$F^2_a= F^2_m\cap CH^p(X)_{\rm alg}$) and $\gr^0_{F_m}CH^p(X)$ has finite
rank.

For $p=1$ this extension is the usual extension associated to the connected
component of the identity of the Picard functor, \ie
$A^1_{X/k}=\Pic_{X/k}^0$ and $\gr^0_{F_m}$ is the N\'eron-Severi of $X$.
If
$p=\dim X$ then $F_m^1$ will be the kernel of the degree map and $F_m^2$
is the
Albanese kernel, \ie $A^{\dim X}_{X/k}$ is the Albanese variety and
$\gr^0_{F_m}= \Z^{\oplus c}$ where $c$ is the number of components of $X$.

However, if $1<p<\dim X$ then $CH^p(X)_{\rm alg}\neq CH^p(X)_{\rm hom}$ in
general.
Let $Grif^p(X)$ denote the quotient group, \ie the Griffiths group of $X$.
Since $\gr^1_{F_a}CH^p(X) = A^p_{X/k}(k)$ and $\gr^0_{F_a}CH^p(X)=NS^p(X)$
we then have a diagram with exact rows and columns
\begin{equation}\label{abelext}
\begin{array}{ccccccc}
&& &0&&0&\\
&& &\downarrow&&\downarrow&\\
&0&\to &F^2_m/F^2_a&\to&Grif^p(X)&\\
&\downarrow& &\downarrow&&\downarrow&\\
0\to &A^p_{X/k}(k)&\to& CH^p(X)/F_a^2& \to & NS^p(X)& \to 0\\
&\downarrow& &\downarrow&&\downarrow&\\
0\to &\gr^1_{F_m}CH^p(X)& \to& CH^p(X)/F_m^2 &\to & \gr^0_{F_m}CH^p(X)&
\to 0\\
&& &\downarrow&&\downarrow&\\
&& &0&&0&
\end{array}
\end{equation}
Note that these extensions still fail to be of the same kind of the Pic
extension.  However, considering the extension $$0\to CH^p(X)_{\rm alg} \to
CH^p(X) \to NS^p(X) \to 0$$ we may hope for a natural regular
homomorphism to $k$-points of an abstract extension in the category of
group schemes (locally of finite type over $k$) $$ 0\to A \to G \to N \to
0$$ where $G$ is a commutative group scheme which is an extension of a
discrete group of finite rank $N$, associated to the abelian group of
codimension p cycles modulo numerical equivalence, by an abelian variety
$A$, isogenous to the universal regular quotient (\cf \cite{GR}). If
$k=\C$ it is easy to see that such extension $G(\C)$ exists trascendentally.

\section{Hodge 1-motives}
Let $k$ be a field, for simplicity, algebraically closed of characteristic
zero. Consider the $\Q$-linear abelian category ${\rm 1-Mot}_k$ of
1-motives
over $k$ with rational coefficients (see \cite{D} and \cite{BSAP}).
Denote $M_{\sQ}$ the isogeny class of a 1-motive $M=[L\to G]$. The
category ${\rm 1-Mot}_k$ contains (as fully faithful abelian
sub-categories) the tensor category of finite dimensional $\Q$-vector
spaces as well as the semi-simple abelian category of isogeny classes of
abelian varieties.
The Hodge realization (see \cite{D} and \cite{BSAP}) is a fully faithful
functor $$T_{\rm Hodge}: {\rm 1-Mot}_{\sC}\into {\rm MHS}\hspace*{1cm}
M_{\sQ}\mapsto T_{\rm Hodge}(M_{\sQ})$$ defining an equivalence of
categories between ${\rm 1-Mot}_{\sC}$ and the abelian sub-category of
mixed $\Q$-Hodge structures of type $\{(0,0), (0,-1), (-1,0), (-1,-1)\}$
such that $\gr_{-1}^W$ is polarizable.  Under this equivalence a
$\Q$-vector space corresponds to a $\Q$-Hodge structure purely of type
$(0,0)$ and an isogeny class of an abelian variety corresponds to a
polarizable $\Q$-Hodge structure purely of type $\{(0,-1), (-1,0)\}$.

\subsection{Algebraic construction}
Let $X$ be a proper scheme over $k$. We perform such a construction for
simplicial
schemes $\Xs$ coming from universal cohomological descent morphisms $\pi :
\Xs\to X$
(\cf \cite{D}, \cite{GN}, \cite{GS}, \cite{HA} and \cite{DJ}).

Let $\Xs$ be such a proper smooth simplicial scheme over the base field
$k$.  By functoriality, the filtration $F_m^jCH^p$ on each component $X_i$
of $\Xs$ yields a complex $$(F_m^jCH^p)^{ \bullet}: \cdots\to
F_m^jCH^p(X_{i-1})\by{\delta^{*}_{i-1}}
F_m^jCH^p(X_{i})\by{\delta^{*}_{i}}F_m^jCH^p(X_{i+1})\to\cdots$$ where
$\delta^{*}_{i}$ is the alternating sum of the pullback along the face maps
$\partial^{k}_{i}: X_{i+1}\to X_i$ for $0\leq k \leq i+1$.

The complex of Chow groups $(CH^p)^{ \bullet}$, induced from the
simplicial structure as above, is filtered by sub-complexes:
$$0\subseteq (F_m^pCH^p)^{ \bullet}\subseteq \cdots \subseteq (F_m^1CH^p)^{
\bullet} \subseteq (F_m^0CH^p)^{ \bullet}= (CH^p)^{ \bullet}.$$
Define $F_m^*(CH^p)^{ \bullet}\df (F_m^*CH^p)^{ \bullet}.$

The extension (\ref{cyclext}) given by the filtration $F_m^*CH^p(X_i)$
on each component $X_i$ of the simplicial scheme $\Xs$, for a fixed $p\geq
0$, yields the following short exact sequence of complexes
\begin{equation}\label{simpext}
0\to \gr^1_{F_m}(CH^p)^{ \bullet} \to (CH^p)^{ \bullet}/F_m^2 \to
\gr^0_{F_m}(CH^p)^{ \bullet} \to 0
\end{equation}
Note that $ \gr^1_{F_m}(CH^p)^{i}$ contains the group of $k$-points of the
abelian variety $A_{X_i/k}^p$ and, moreover $\gr^0_{F_m}(CH^p)^{i}_{\sQ}$
is the
finite dimensional vector space of codimension p cycles on $X_i$ modulo
homological (or numerical) equivalence.

 From (\ref{simpext}) we then get a long exact sequence of homology groups
and, in particular, we obtain boundary maps $$\lambda^i_m
:H^{i}(\gr^0_{F_m}(CH^p)^{\bullet}) \to
H^{i+1}(\gr^1_{F_m}(CH^p)^{\bullet}).$$ Denote $A_{\Xs/k}^p$ the complex of
abelian varieties $A_{X_i/k}^p$.  Since $A_{\Xs/k}^p(k)$ is a sub-complex
of $\gr^1_{F_m}(CH^p)^{\bullet}$ we then get induced (functorial) maps on
homology groups $$\theta^i: H^{i+1}(A_{\Xs/k}^p(k))\to
H^{i+1}(\gr^1_{F_m}(CH^p)^{\bullet}).$$ Note that (\ref{simpext}) is
involved in a functorial diagram (\ref{abelext}). The corresponding
complex of N\'eron-Severi groups $(NS^p)^{\bullet}$ yield boundary maps
\begin{equation}\label{bound}
\lambda^i_a :H^{i}((NS^p)^{\bullet}) \to H^{i+1}(A_{\Xs/k}^p(k)).
\end{equation}
These maps fit into the following
commutative square $$\begin{array}{ccc} H^{i}(\gr^0_{F_m}(CH^p)^{\bullet})&
\longby{\lambda^i_m}& H^{i+1}(\gr^1_{F_m}(CH^p)^{\bullet})\\
\uparrow& &\uparrow{\ \scriptsize \theta^i}\\
H^{i}((NS^p)^{\bullet})&\longby{\lambda^i_a}&H^{i+1}(A_{\Xs/k}^p(k)).
\end{array}$$
Moreover, the kernel of $\theta^i$ is clearly equal to the image of the
boundary map $$\tau^i :H^{i}(\gr^1_{F_m}(CH^p)^{\bullet}/A_{\Xs/k}^p(k))\to
H^{i+1}(A_{\Xs/k}^p(k)).$$
\begin{question} Is the image of the connected component of the identity of
$H^{i+1}(A_{\Xs/k}^p)$ under $\theta^i$ an abelian variety,
\eg is $\tau^i = 0$ up to a finite group ?
\end{question}
This is clearly the case if $p=1$ (see below for the case $k=\C$).
For $k=\C$ this question is related to Griffiths Problem~E in
\cite{GR} asking a description of the ``invertible points'' of the
intermediate jacobians (also \cf Mumford-Griffiths Problem~F in
\cite{GR}).

Assume that the above question has a positive answer and denote
$H^{i+1}(A_{\Xs/k}^p)^{\dag}$ the so obtained abelian variety.  We then
obtain an algebraically defined 1-motive as follows.

Let $H^{i}(\gr^0_{F_m}(CH^p)^{\bullet})^{\dag}$ be the sub-group of those
elements in $H^{i}(\gr^0_{F_m}(CH^p)^{\bullet})$ which are mapped to
$H^{i+1}(A_{\Xs/k}^p)^{\dag}$ under the boundary map $\lambda^i_m$
above.
\begin{defn} {\rm Let $\Xs$ be such a smooth proper simplicial scheme over
$k$.
Denote
$$ \Xi^{i,p}\df [ H^{i}(\gr^0_{F_m}(CH^p)^{\bullet})^{\dag}
\longby{\lambda^i_m}
H^{i+1}(A_{\Xs/k}^p)^{\dag}]_{\sQ}$$
the isogeny 1-motive obtained from the construction above.
Call $ \Xi^{i,p}$ the {\em Hodge 1-motive} of the simplicial scheme.}
\end{defn}
We expect that $\Xi^{i,p}$ is independent of the choice of $\pi : \Xs\to
X$. The main motivation for questioning the existence of such a purely
algebraic construction is given by the following analytic counterpart.

\subsection{Analytic construction}
Let ${\rm MHS}$ be the abelian category of usual Deligne's mixed
Hodge structures \cite{D}. An object $H$ of ${\rm MHS}$ is defined as a
triple $H = (H_{\sZ}, W, F)$ where $H_{\sZ}$ is a finitely generated
$\Z$-module, $W$ is a finite increasing filtration on $H_{\sZ}\otimes\Q$
and $F$ is a
finite decreasing filtration on $H_{\sZ}\otimes\C$ such that $W, F$ and
$\bar F$
is a system of opposed filtrations.

Let $H\in {\rm MHS}$ be a torsion free mixed Hodge structure with positive
weights.  Let $W_*H$ denote the sub-structures defined by the intersections
of the weight filtration and $H_{\sZ}$.  Let $p$ be a fixed
positive integer and assume that $\gr^W_{2p-1}H$ is polarizable.

Consider the following extension in the abelian category ${\rm MHS}$
\begin{equation}\label{ext} 0\to \gr^W_{2p-1}H
\to \frac{W_{2p}H}{W_{2p-2}H}\to \gr^W_{2p}H \to 0 \end{equation} Taking
$\Hom (\Z(-p), -)$ we get the extension class map $$e^p: \Hom
(\Z(-p),\gr^W_{2p}H)\to \Ext (\Z(-p), \gr^W_{2p-1}H)$$ where $\Hom
(\Z(-p),\gr^W_{2p}H)=H_{\sZ}^{p,p}$ is the sub-structure of $(p,p)$-classes
in $\gr^W_{2p}H$ and $$ \Ext (\Z(-p), \gr^W_{2p-1}H)\cong J^p(H)\df
\frac{\gr^W_{2p-1}H_{\sC}}{F^p +\gr^W_{2p-1}H_{\sZ}}$$ is a compact complex
torus.  Note that (\cf \cite{CA}) $$\Ext (H_{\sZ}^{p,p},
\gr^W_{2p-1}H)\cong \Hom
(H_{\sZ}^{p,p},J^p(H)).$$ Thus $e^p\in \Hom (H_{\sZ}^{p,p}, J^p(H))$
corresponds to a unique extension class \begin{equation}\label{he} 0\to
\gr^W_{2p-1}H \to H^e \to H_{\sZ}^{p,p}\to 0 \end{equation} which is the
pull-back extension associated to $H_{\sZ}^{p,p}\into \gr^W_{2p}H$ and
(\ref{ext}).  Moreover, since we always have $\gr^W_{2p-1}H\cap F^p =0$
then $$F^p\cap H^e_{\sZ} = \ker (H_{\sZ}^{p,p}\by{e^p}J^p(H)).$$ Now, if
$\gr^W_{2p-1}H$ is (polarizable) of level $1$ then the torus $J^p(H)$ is an
abelian variety and $H^e$ is the Hodge realization of the 1-motive over
$\C$ defined by the extension class map $e^p$ above.

In general, let $H^{\prime}$ be the largest sub-structure of $W_{2p-1}H$
which is polarizable and purely of type $\{(p-1,p), (p,p-1)\}$ modulo
$W_{2p-2}H$, \ie if $$H^{2p-1}_a \df (H^{p-1,p}+ H^{p,p-1})_{\sZ}$$ is the
polarizable sub-structure of $\gr^W_{2p-1}H$ of those elements which are purely of the
above type, then $H^{\prime}$ is defined by the following pull-back
extension $$0\to W_{2p-2}H \to H^{\prime}\to H^{2p-1}_a\to 0,$$ along the
canonical projection $W_{2p-1}H\onto \gr^W_{2p-1}H.$

Let $H^{\prime\prime}\subseteq W_{2p}H$ be defined by the following
pull-back extension $$0\to W_{2p-1}H \to H^{\prime\prime}\to H_{\sZ}^{p,p}
\to 0,$$  along the canonical projection $W_{2p}H\onto \gr^W_{2p}H.$
Thus, the extension (\ref{he}) can be regarded as the push-out of such
extension involving $H^{\prime\prime}$ along  $W_{2p-1}H\onto \gr^W_{2p-1}H.$
Namely, we obtain that
$$\frac{H^{\prime\prime}}{W_{2p-2}H}=H^e$$ fitting in the extension
\begin{equation}\label{heprime} 0\to \frac{H^{\prime}}{W_{2p-2}H} \to H^e
\to \frac{H^{\prime\prime}}{H^{\prime}}\to 0.  \end{equation}
Let $$h^p: \Hom (\Z(-p),\frac{H^{\prime\prime}}{H^{\prime}})\to \Ext (\Z(-p),
\frac{H^{\prime}}{W_{2p-2}H})$$ be the extension class map obtained from
(\ref{heprime}).  \begin{propose} The map $h^p$ above yields a 1-motive
over $\C$ which is just the restriction of $e^p:H_{\sZ}^{p,p}\to J^p(H)$ to
the largest abelian subvariety in $J^p(H)$.  In particular: $\ker (h^p) =
\ker (e^p).$ \end{propose} \begin{proof} Since $H^{\prime}/W_{2p-2}H=
H_a^{2p-1}$ by construction we have that $ \Ext
(\Z(-p),H_a^{2p-1})$ is the largest abelian
sub-variety of $ \Ext (\Z(-p),\gr^W_{2p-1}H )=J^p(H)$.  Moreover, note that
we also have induced extensions $$0\to \frac{H^{\prime}}{W_{2p-2}H} \to
\gr^W_{2p-1}H \to \frac{W_{2p-1}H}{H^{\prime}}\to 0$$ and $$0\to
\frac{W_{2p-1}H}{H^{\prime}} \to \frac{H^{\prime\prime}}{H^{\prime}}\to
H_{\sZ}^{p,p}\to 0$$ yielding, together with (\ref{he}) and
(\ref{heprime}), the following commutative diagram with exact rows
$$\begin{array}{c} \Hom (\Z(-p),\frac{W_{2p-1}H}{H^{\prime}})\to \Hom
(\Z(-p),\frac{H^{\prime\prime}}{H^{\prime}})\to \Hom
(\Z(-p),H_{\sZ}^{p,p})\to \Ext (\Z(-p),\frac{W_{2p-1}H}{H^{\prime}})\\ \veq
\hspace*{3cm}h^p\downarrow \hspace*{3cm}\downarrow e^p\hspace*{3cm} \veq\\
\Hom (\Z(-p),\frac{W_{2p-1}H}{H^{\prime}})\to \Ext (\Z(-p),
\frac{H^{\prime}}{W_{2p-2}H})\to\Ext (\Z(-p), \gr^W_{2p-1}H)\to\Ext
(\Z(-p),\frac{W_{2p-1}H}{H^{\prime}}) \end{array}$$
Since $\Hom (\Z(-p),\frac{W_{2p-1}H}{H^{\prime}}) = 0 = \Hom (\Z(-p),H_a^{2p-1})$
everything follows from a diagram chase.  \end{proof}
\begin{defn}{\rm Let $A^p(H)\df \Ext
(\Z(-p),H_a^{2p-1})$ denote the abelian part of the
compact complex torus $J^p(H)$.  Denote $H^p(H)\df \Hom
(\Z(-p),H^{\prime\prime}/H^{\prime})$ the group of Hodge cycles, \ie the
sub-group of $H_{\sZ}^{p,p}$ mapped to $A^p(H)$ under the extension class
map $e^p$.  Define $$H^h\df T_{\rm Hodge}([H^p(H)\by{h^p}A^p(H)])$$ the
mixed Hodge structure corresponding to the 1-motive defined from
(\ref{heprime}) above. Call this 1-motive the {\it Hodge 1-motive} of $H$.}
\end{defn}
We remark that the mixed Hodge structure $H^h$ clearly corresponds to the
following extension
$$0\to H_a^{2p-1}\to H^h\to H^p(H)\to 0,$$
obtained by pulling back $H^p(H) = F^p\cap
(H^{\prime\prime}/H^{\prime})_{\sZ}$ along the projection $H^e\onto
H^{\prime\prime}/H^{\prime}$ in (\ref{heprime}). In particular
$H^h\subseteq H^e$ and
$F^p\cap H_{\sZ}\subseteq F^p\cap H^h_{\sZ} = F^p\cap H^e_{\sZ}.$

\subsection{Hodge conjecture for singular varieties}
Let $X$ be a proper smooth $\C$-scheme.
The coniveau or arithmetic filtration (\cf \cite{GH})
$$N^iH^j(X) \df \ker (H^j(X) \to  \limdir{{\rm codim}_XZ\geq i} H^j
(X-Z))$$
yields a filtration by Hodge sub-structures of $H^j(X)$. We have that
$N^iH^j(X)$ is of level $j-2i$ and
\begin{equation}\label{coin}
N^iH^j(X)_{\sQ}\subseteq H^j(X, \Q)\cap F^iH^j(X).
\end{equation}
\begin{conj} {\em (Grothendieck-Hodge conjecture \cite{GH})}
The left hand side of (\ref{coin}) is the largest sub-space of the right
hand side, generating a sub-space of $H^j(X, \C)$ which is a sub-Hodge
structure.
\end{conj}

Let $X$ be a proper (integral) $\C$-scheme.  Recall that the weight
filtration on $H^{*}(X, \Q)$ is given by the canonical spectral sequence of
mixed $\Q$-Hodge structures $$E_1^{s, t}=H^{t}(X_s)\implies
\HH^{s+t}(\Xs)$$ for any smooth and proper hypercovering $\pi : \Xs \to X$
and universal cohomological descent $H^{*}(X)\cong \HH^{*}(\Xs)$ (see
\cite{D}).  In fact, the spectral sequence degenerates at $E_2$.  Denote
$(H^{t})^{\bullet}$ the complexes $E_1^{\cdot, t}$ of $E_1$-terms.  We
clearly have $$H^i((H^{t})^{\bullet})=\gr^W_{t}H^{t+i}(X).$$ Consider the
complexes $(N^lH^{t})^{\bullet}$ induced from the coniveau filtration
$N^lH^{t}(X_{i})$ on the components $X_i$ of the simplicial scheme $\Xs$.
We then have a natural map of Hodge structures $$ \nu^{i,l}
:H^i((N^lH^{t})^{\bullet}) \to \gr^W_{t}H^{t+i}(X).$$ Note that the image
of $\nu^{i,l}$ is contained in the sub-space $\gr^W_{t}H^{t+i}(X, \Q)\cap F^l$.
\begin{conj}\label{GHC}
The image of $\nu^{i,l}$ is the largest sub-space of
$\gr^W_{t}H^{t+i}(X,\Q)\cap F^l$ which is a sub-Hodge structure of
$\gr^W_{t}H^{t+i}(X)$.
\end{conj}
It is reasonable to expect that such a statement will follow from the
original Grothendieck-Hodge conjecture and abstract Hodge theory.

Grothendieck-Hodge conjecture (for coniveau $p$ and degrees $2p, 2p+1$) can
be reformulated as follows (\cf Grothendieck's remark on
motives in \cite{GH}).  Consider $\gr^0_{F_m}CH^{p}(X)$ and
$A^{p+1}_{X/k}\subseteq \gr^1_{F_m}CH^{p+1}(X)$ (for $k=\C$ this is the
algebraic part of $J^{p+1}(X)$) as 1-motives with rational coefficients.
The Hodge realization of these algebraically defined 1-motives is
$N^pH^{2p}(X)$ and $N^pH^{2p+1}(X)$ respectively.
\begin{conj} Let $X$ be smooth and proper over $\C$. Then
$$ T_{\rm Hodge}([\gr^0_{F_m}CH^{p}(X)\to 0]_{\sQ}) = H^{p,p}_{\sQ}$$
and $$T_{\rm Hodge}([0\to
A^{p+1}_{X/k}]_{\sQ})=(H^{p,p+1}+H^{p+1,p})_{\sQ}.$$
\end{conj}
Note that $\gr^0_{F_m}CH^{p}(X)$ would be better defined as
$\cZ^p(X)/\equiv_{\rm num}$, up to torsion.

Now apply to the mixed $\Q$-Hodge structure $H=H^{2p+i}(X)$ the
construction
performed in the previous section. For a fixed pair $(i,p)$ of integers
recall that (\ref{ext}) is an extension of $\gr^W_{2p}H^{2p+i}(X)$
by $\gr^W_{2p-1}H^{2p+i}(X)$, where:
$$H^{i+1}((H^{2p-1})^{\bullet})=\gr^W_{2p-1}H^{2p+i}(X)=
\frac{\ker (H^{2p-1}(X_{i+1}) \rightarrow H^{2p-1}(X_{i+2}))}{\im
(H^{2p-1}(X_{i}) \rightarrow H^{2p-1}(X_{i+1}))}$$ and
$$H^{i}((H^{2p})^{\bullet})=\gr^W_{2p}H^{2p+i}(X)=
\frac{\ker(H^{2p}(X_{i}) \rightarrow H^{2p}(X_{i+1}))}{\im
(H^{2p}(X_{i-1}) \rightarrow H^{2p}(X_{i}))}.$$
We then have that $J^p(H)= J^p(H^{i+1}((H^{2p-1})^{\bullet})))$ is
isogenous to
$H^{i+1}((J^p)^{\bullet}))$ where $(J^p)^{\bullet}$ is the complex of
jacobians $J^p(X_i)$ of the components $X_i$.

Consider the complex $A_{\Xs/\sC}^p$ of abelian sub-varieties given by the
algebraic part of intermediate jacobians.
The complex $A_{\Xs/\sC}^p(\C)$ is a sub-complex of the complex of
compact tori $(J^p)^{\bullet}$. Therefore, the induced maps
$$H^{i+1}(A_{\Xs/\sC}^p(\C)) \to J^p(H)$$
are holomorphic mappings (which factor through $\theta^i$) and whose image
is isogenous to an abelian sub-variety of the maximal abelian sub-variety
$A^p(H)$ of
$J^p(H)$.

Moreover, the homology of the complex of $(p,p)$-classes is mapped to
$H^{p,p}_{\sQ} = F^p\cap H^{i}((H^{2p})^{\bullet})$. Thus, there
are canonical maps
$$H^{i}((NS^p)^{\bullet}) \to H^{p,p}_{\sQ}$$
which factor through $H^{i}(\gr_{F_m}^0(CH^p)^{\bullet})$.

We expect that the Hodge 1-motive of the simplicial
scheme $\Xs$ would be canonically isomorphic to the Hodge 1-motive of
$H^{2p+i}(X)$.
\begin{conj} \label{MHC}
Let $X$ be a proper $\C$-scheme and let $\pi : \Xs \to X$ be a smooth and
proper hypercovering. Let $H^{2p+i}(X)$ denote Deligne's mixed $\Q$-Hodge
structure
on $H^{2p+i}(X, \Q)$, \ie obtained by the universal cohomological descent
isomorphism $H^{2p+i}(X, \Q)\cong \HH^{2p+i}(\Xs , \Q)$.
\begin{enumerate}
\item The following square
$$\begin{array}{ccc}
H^{i}((NS^p)^{\bullet})&\by{\lambda^i_a}& H^{i+1}(A_{\Xs/\sC}^p(\C))\\
\downarrow & &\downarrow\\
H^{2p+i}(X)^{p,p} & \by{e^p} & J^p(H^{2p+i}(X))
\end{array}$$
commutes, yielding a motivic ``cycle class map''.
\item The image of the motivic ``cycle class map'' is the Hodge
1-motive of the $\Q$-Hodge structure $H^{2p+i}(X)$.
\item We have that
$$T_{\rm Hodge}(\Xi^{i,p}) \cong H^{2p+i}(X)^h.$$
\end{enumerate}
\end{conj}

\begin{remark}{\em Note that if $X$ is smooth and proper then $H^{2p+i}(X)$
is pure and $H^{2p+i}(X)^h\neq 0$ if and only if $i= 0, -1$ ($p$ fixed).
In this case, the above conjecture follows from the reformulation of
Grothendieck-Hodge conjecture for $H^{2p}(X)$ and $H^{2p-1}(X)$.}
\end{remark}

\section{Local Hodge theory}
See \cite{D} for notations, definitions and properties of mixed Hodge
structures.

\subsection{Infinite dimensional mixed Hodge structures}
Let ${\rm MHS}$ denote the abelian category of usual Deligne's $A$-mixed
Hodge structures \cite{D}, \ie for $A$ a noetherian subring of $\R$ such
that $A\otimes\Q$ is a field, an object $H$ of ${\rm MHS}$ is defined as a
triple $H = (H_A, W, F)$ where $H_A$ is a finitely generated $A$-module,
$W$ is a finite increasing filtration on $H_A\otimes\Q$ and $F$ is a finite
decreasing filtration on $H_A\otimes\C$ such that $W, F$ and $\bar F$ is a
system of opposed filtrations.

\begin{defn}{\em An $\infty$-mixed Hodge structure $H$ is a triple $(H_A,
W, F)$ where $H_A$ is any $A$-module, $W$ is a finite increasing
filtration on
$H_A\otimes\Q$ and $F$ is a finite decreasing filtration on $H_A\otimes\C$
such that $W, F$ and $\bar F$ is a system of opposed filtrations.

Denote ${\rm MHS}^{\infty}$ the category of $\infty$-mixed Hodge
structures or ``infinite dimensional'' mixed Hodge structures, \ie where
the morphisms are those which are compatible with the
filtrations.}\end{defn}

The category ${\rm MHS}^{\infty}$ is abelian and ${\rm MHS}$ is a fully
faithful abelian subcategory of ${\rm MHS}^{\infty}$.
Note that similar categories of infinite dimensional mixed Hodge structures
already appeared in the literature, see Hain \cite{HN} and Morgan
\cite{MO}.  For example the category of limit mixed Hodge structures ${\rm
MHS}^{\lim}$, \ie whose objects and morphisms are obtained by formally add
to ${\rm MHS}$ (small) filtered colimits of objects in ${\rm MHS}$ with
colimit morphisms.

Consider the case $A=\Q$.  In this case, in the category ${\rm
MHS}^{\infty}$ we have infinite products of those families of objects
$\{H_i\}_{i\in I}$ such that the induced families of filtrations
$\{W_{i}\}_{i\in I}$ and $\{F_{i}\}_{i\in I}$ are finite.
Moreover such a (small) product of epimorphisms is an epimorphism.

For the sake of exposition we often call mixed Hodge structures the
objects of ${\rm MHS}$ as well as those of ${\rm  MHS}^{\infty}$ (or ${\rm
MHS}^{\lim}$).

\subsection{Zariski sheaves of mixed Hodge structures}
Let $X$ denote the (big or small) Zariski site on an algebraic $\C$-scheme.
However, most of the results in this section are available for any
topological space or Grothendieck site.

Denote $X_{\d}$ a simplicial object of the category of algebraic
$\C$-schemes over $X$: recall that (see \cite[5.1.8]{D}) simplicial sheaves
on $X_{\d}$ can be regarded as objects of a Grothendieck topos with enough
points.

Consider presheaves of mixed Hodge structures. Note that a presheaf of
usual
Deligne's $A$-mixed Hodge structures will have its stalks in ${\rm
MHS}^{\lim}$.

Consider those presheaves (resp. simplicial presheaves) of $\Q$-mixed Hodge
structures on $X$ (resp. on $X_{\d}$) such that the filtrations are finite
as filtrations of sub-presheaves on $X$. These presheaves can be sheafified
to sheaves having finite filtrations and preserving the above conditions on
the
stalks.

Make the following working definition of sheaves (or simplicial
sheaves) of mixed Hodge structures. Let $A, \Q$ and $\C$ denote as well the
constant sheaves on $X$ (or $X_{\d}$) associated to the ring $A$, the
rationals and the complex numbers.

\begin{defn}{\rm A (simplicial) sheaf $\cH$ of $A$-mixed Hodge structures,
or
``$A$-mixed sheaf'' for short, is given by the following set of datas:
\begin{description}
\item[{\it i)}] a (simplicial) sheaf $\cH_A$ of $A$-modules,
\item[{\it ii)}] a finite (exhaustive) increasing filtration $\cW$ by
$\Q$-subsheaves
of $\cH_{\sQ}\df \cH_A\otimes\Q$,
\item[{\it iii)}] a finite (exhaustive) decreasing filtration $\cF$ by
$\C$-subsheaves
of $\cH_{\sC}\df \cH_A\otimes\C$;
\end{description}
satisfying the condition that $\cW,\cF$ and $\bar\cF$ is a system of
opposed filtrations, \ie we have that $$gr^p_{\cF}gr^q_{\bar
\cF}gr_{n}^{\cW}(\cH) =  0$$ for $p+q\neq n$.}
\end{defn}

There is a canonical decomposition  $$gr_n^{\cW}(\cH) =
\bigoplus_{p+q=n}^{}
\cA^{p,q}$$ where $\cA^{p,q} \df \cF^p\cap \bar \cF^q$ and conversely.

In the case of a simplicial sheaf assume that the filtrations are given by
simplicial subsheaves, \ie the simplicial structure should be compatible
with the filtrations on the components.  A simplicial $A$-mixed sheaf
$\cH_{\Xs}$ on the simplicial space $\Xs$ can be regarded (\cf
\cite[5.1.6]{D})
as a family of $A$-mixed sheaves $\cH_{X_i}$ (on the components $X_i$) such
that the simplicial structure is also compatible with the filtrations
$\cW_{X_i}$
and $\cF_{X_i}$ of $\cH_{X_i}$.

A morphism of $A$-mixed sheaves is a morphism of (simplicial) sheaves of
$A$-modules which is compatible with the filtrations. Denote
${\cal  MHS}_X$ and ${\cal MHS}_{X_{\d}}$ the corresponding categories.

In order to show the following Lemma one can just use Deligne's Theorem
\cite[1.2.10]{D}.
\begin{lemma}\label{abel} The categories ${\cal MHS}_X$ and ${\cal
MHS}_{X_{\d}}$ of $\Q$-mixed sheaves are abelian categories. The
kernel (resp. the cokernel) of a morphism $\varphi : \cH \to \cH'$ has
underlying  $\Q$  and $\C$-vector spaces the kernels (resp. the cokernels)
of $\varphi_{\sQ}$ and $\varphi_{\sC}$ with induced filtrations. Any
morphism is strictly compatible with the filtrations. The functors
$gr_{\cW}$ and  $gr_{\cF}$ are exacts.
\end{lemma}

Note that if $X=\{\infty\}$ is the singleton then ${\cal
MHS}_{\infty}$ is equal to ${\rm  MHS}^{\infty}$.
Examples of $\Q$-mixed sheaves are clearly given by constant sheaves
associated
to  $\Q$-mixed Hodge structures, yielding a canonical fully faithful
functor
$${\rm MHS}^{\infty}\to {\cal MHS}_X.$$

Stalks of a $\Q$-mixed sheaf $\cH$ are in ${\rm MHS}^{\infty}$, the
filtrations being induced stalkwise. In fact, the condition on the
filtrations given with any  $\Q$-mixed sheaf is local, at any point of $X$.
Skyscraper sheaves $x_*(H)$ associated to an object $H\in{\rm
MHS}^{\infty}$
and a point $x$ of $X$ are in ${\cal MHS}_X$. There is a natural
isomorphism
$$\Hom_{{\rm MHS}^{\infty}}(\cH_x, H) \cong \Hom_{{\cal MHS}_X}(\cH,
x_*(H)).$$

More generally, a presheaf in ${\rm MHS}^{\infty}$, with finite filtrations
presheaves, can be sheafified to an $A$-mixed sheaf, in a canonical way, by
applying the usual sheafification process to the filtrations together with
the presheaf.

\begin{defn}{\rm Say that an $A$-mixed sheaf $\cH$ is flasque if $\cH_A$
is a
flasque sheaf.} \end{defn}

For a given $\Q$-mixed sheaf $\cH$ we then dispose of a
canonical flasque $\Q$-mixed sheaf $$\prod_{x\in X}x_*(\cH_x)$$
where the product is taken over a set of points of $X$.

\subsection{Hodge structures and Zariski cohomology}
We show that, if $\cH$ is a $\Q$-mixed sheaf then there is a unique
$\Q$-mixed Hodge
structure on the sections such that $\Gamma (-, gr(\dag)) = \gr\Gamma
(-,\dag)$.
In  fact, the mixed Hodge structure on the $\Q$-vector space of (global)
sections is such  that the following $$\Gamma (X,\cH_{\sQ}) \subset
\prod_{x \in
X}\cH_x$$ is strictly compatible with the filtrations; in the same way,
for a simplicial
sheaf, the following
$$\Gamma (X_{\d},\cH_{\sQ,_{\d}}) \subset \ker \prod_{x \in X_0}\cH_x\to
\prod_{x \in X_1}\cH_x$$
is strictly compatible with the filtrations.

\begin{propose} Let $\cH_{X}\in{\cal MHS}_X$ and $\cH_{X_{\d}}\in{\cal
MHS}_{X_{\d}}$ as above. There are left exact functors
$$\cH_{X}\mapsto \Gamma (X,\cH_X) \mbox{\hspace*{1cm}} {\cal MHS}_X \to
{\rm  MHS}^{\infty}$$ and
$$\cH_{X_{\d}}\mapsto\Gamma (X_{\d},\cH_{X_{\d}}) \mbox{\hspace*{1cm}}
{\cal  MHS}_{X_{\d}} \to {\rm MHS}^{\infty}$$
These functors yield $\Q$-mixed Hodge structures on the usual cohomology,
which we
denote by $H^*(X,\cH_X)$ and $H^*(X_{\d},\cH_{X_{\d}})$ respectively, such
that if
$$0\to \cH' \to \cH\to \cH'' \to 0$$ is exact in ${\cal MHS}_X$,
respectively in ${\cal  MHS}_{X_{\d}}$, then
$$\cdots\to H^i(X,\cH_X)\to H^i(X,\cH''_X)\to H^{i+1}(X,\cH'_X)\to \cdots$$
is exact in ${\rm MHS}^{\infty}$, and respectively for the cohomology of
$X_{\d}$: moreover, in the latter case we have a spectral sequence
$$E^{p,q}_1 = H^q(X_p,\cH_{X_p}) \implies \HH^{p+q}(X_{\d},\cH_{X_{\d}})$$
in the category ${\rm MHS}^{\infty}$.
\end{propose}
\begin{proof} In fact, there is an extension in ${\cal MHS}$
$$0\to \cH \to \prod_{x\in X}x_*(\cH_x)\to \cZ^1\to 0$$
where $\cZ^1$ has the quotient $\Q$-mixed structure; as usual, we then get
another extension
$$0\to \cZ^1\to \prod_{x\in X}x_*(\cZ^1_x)\to \cZ^{2}\to 0$$
and so on. We therefore get a flasque resolution
$$\prod_{x\in X}x_*(\cH_x)\to \prod_{x\in X}x_*(\cZ^1_x)\to
\prod_{x\in X}x_*(\cZ^2_x)\to \cdots$$
in ${\cal MHS}$: the canonical $\Q$-mixed flasque resolution.

If $\cH_X$ is flasque then $\Gamma (X,\cH_X)$ has a canonical
$\Q$-mixed Hodge structure as claimed above; note that the filtrations
would be given by flasque sub-sheaves.

In general, by construction, the cohomology is the homology of the complex
of
sections in ${\rm MHS}^{\infty}$.  Thus $H^*(X,\cH_X)$ has a canonical
$\Q$-mixed Hodge structure.  The same argument applies to the
total complex of the double complex of flasque $\Q$-mixed sheaves given by
the canonical resolutions on each component of a simplicial sheaf.

Refer to [SGA4] and \cite[Chapter IV]{LM} for a
construction of  canonical Godement resolutions available on any site and
compare  \cite[1.4.11]{D} for the existence of bifiltered resolutions.
\end{proof}

In particular, the mixed Hodge structure $H^* (X,\cH_X)$ is such that $H^*
(X,\cH_X)_{\sQ} = H^* (X,\cH_{\sQ})$, $WH^* (X,\cH_X)_{\sQ} = H^*
(X,\cW\cH_{\sQ})$ and $FH^* (X,\cH_X)_{\sC} = H^* (X,\cF\cH_{\sC})$.  There
is a decomposition $$\gr_n^{W}H^* (X, \cH_X) = \bigoplus_{p+q=n}^{} H^* (X,
\cA^{p,q}_X).$$

\begin{remark}{\em Note that any (non-canonical) $\Q$-mixed flasque
resolution in ${\cal MHS}$ yields a bifiltered complex and a bifiltered
quasi-isomorphism with the canonical resolution. Therefore, the so obtained
${\infty}$-mixed Hodge structure on the cohomology is unique up to
isomorphism.}
\end{remark}

Considering complexes in ${\cal MHS}_X$ and ${\cal
MHS}_{X_{\d}}$ we construct the derived categories of $\Q$-mixed sheaves
$\cD^*({\cal MHS}_X)$ and $\cD^*({\cal MHS}_{X_{\d}})$ as usual, as well as
$\cD^*({\rm MHS}^{\infty})$. We have a total derived functor
$$R\Gamma (X, -): \cD^*({\cal MHS}_X)\to \cD^*({\rm MHS}^{\infty})$$
sending a complex of $\Q$-mixed sheaves to the total complex of sections of
its canonical resolution. Moreover, if $f:X\to Y$ is a continuous map, we
have a higher direct image $\Q$-mixed sheaf $R^qf_*(\cH_X)$ on $Y$, which
corresponds as well to an exact functor
$$Rf_*:\cD^*({\cal MHS}_X)\to \cD^*({\cal MHS}_Y).$$
There is an inverse image exact functor
$f^*:\cD^*({\cal MHS}_Y)\to \cD^*({\cal MHS}_X).$
Moreover, Grothendieck six standard operations can be otained in the
derived category of $\Q$-mixed sheaves.

\subsection{Local-to-global properties}
Let $X$ be an algebraic $\C$-scheme and let $X_{\rm an}$ be the
associated analytic space.  For any Zariski open subset $U\subseteq X$ the
corresponding integral cohomology $H^r(U_{\rm an},\Z(t))$ carry a mixed
Hodge
structure (see \cite[8.2]{D}) such that the restriction maps $H^r(U_{\rm
an},\Z (t))\to H^r(V_{\rm an},\Z (t))$ for $V\subseteq U$ are maps of
mixed Hodge
structures.  Thus the presheaf of mixed Hodge structures
\begin{equation}\label{sheaf}
U \mapsto H^r(U_{\rm an},\Q (t))
\end{equation} can be sheafified to a Zariski $\Q$-mixed sheaf.  In
fact, for a fixed $r$, the resulting non-zero Hodge numbers of $H^r(U_{\rm
an},\Q)$, for any $U$, are in the finite set $[0,r]\times [0,r]$ (see
\cite[8.2.4]{D}).
\begin{defn} {\em Denote $\cH^r_X(\Q (t))$ the $\Q$-mixed sheaf obtained
hereabove. For $\Xs$ a simplicial $\C$-scheme denote $\cH^r_{\Xs}$ the
simplicial
$\Q$-mixed sheaf given by $\cH^r_{X_p}$ on the component $X_p$. }
\end{defn}

If $X$ has algebraic dimension $n$ then all its Zariski open affines $U$
do have
dimension $\leq n$ thus $\cH^r_X=0$ for $r>n$.

\begin{cor}\label{zarhodge}
The Zariski cohomology groups $H^*(X,\cH^r_X)$ carry
$\infty$-mixed Hodge structures. Possibly non-zero Hodge numbers of
$H^*(X,\cH^r_X)$
are in the finite set $[0,r]\times [0,r]$. The Zariski cohomology
$\HH^{*}(X_{\d},\cH_{X_{\d}}^r)$ carry
$\infty$-mixed Hodge structure and the canonical spectral sequence
$$E^{p,q}_1 = H^q(X_p,\cH_{X_p}^r) \implies
\HH^{p+q}(X_{\d},\cH_{X_{\d}}^r)$$
is in the category ${\rm MHS}^{\infty}$.
\end{cor}

Let $\omega : X_{\rm an} \to X_{\rm Zar}$ be the continuous map of sites
induced by the identity mapping. We then have a Leray spectral sequence
$$L^{q,r}_2 = H^q(X_{\rm Zar},R^r\omega_*(\Z))\implies H^{q+r}(X_{\rm
an},\Z)$$
of abelian groups. Since $R^r\omega_*(\Q)\cong \cH_{X}^r$ these sheaves
can be
regarded as $\Q$-mixed sheaves and its Zariski cohomology carry
$\infty$-mixed
Hodge structures as above.

For $\Xs$ a simplicial scheme we thus have $\omega_{\d} : (X_{\d})_{\rm
an} \to
(X_{\d})_{\rm Zar}$ as above and a Leray spectral sequence
$$L^{q,r}_2 = \HH^q(X_{\d}, R^r(\omega_{\d})_*(\Z))\implies
\HH^{q+r}(X_{\d},\Z)$$
where $R^r(\omega_{\d})_*(\Z)\cong \cH_{\Xs}^r$.

\begin{thm}\label{l2g} {\em (Local-to-global)} There are spectral
sequences
$$L^{q,r}_2 = H^q(X,\cH_{X}^r)\implies H^{q+r}(X_{\rm an},\Q)$$
and $$L^{q,r}_2 = \HH^q(X_{\d}, \cH_{\Xs}^r)\implies
\HH^{q+r}(X_{\d},\Q)$$ in the category of $\infty$-mixed Hodge structures.
\end{thm}

The proof of this compatibility result will appear elsewhere; however, for
smooth $\C$-schemes and using (\ref{ares}) below, the compatibility follows
from \cite[Corollary 4.4]{KP}.

\section{Edge maps}
Recall that the classical cycle class maps can be obtained {\em via}\, edge
homomorphisms in the coniveau spectral sequence. This is a
consequences of Bloch's formula \cite{BO}.  Working simplicially we then
construct certain cycle class maps for singular varieties {\em via}\, edge
maps in the local-to-global spectral sequence.  We first show
that the results of \cite{BO} hold in the category of $\infty$-mixed Hodge
structures.

\subsection{Bloch-Ogus theory}
 From Deligne \cite[8.2.2 and 8.3.8]{D} the cohomology groups
$H^*_Z(X)$ (= $H^*(X {\rm mod} X-Z,\Z)$ in Deligne's notation) carry a
mixed
Hodge structure fitting into long exact sequences
 \begin{equation}\label{loc}
\cdots \to H_Z^j(X)\to H_T^j(X)\to
H_{T-Z}^j(X-Z)\to H_Z^{j+1}(X)\to \cdots
\end{equation} for any pair $Z\subset T$ of closed subschemes of $X$.

Since classical Poincar\'e duality is compatible with the mixed Hodge
structures involved, then the functors
$$Z\subseteq X \mapsto (H_Z^*(X),H_*(Z))$$
yield a Poincar\'e duality theory with supports (see \cite{BO} and
\cite{JA})
in the abelian tensor category of mixed Hodge structures.
Furthermore we have that the above theory is appropriate for algebraic
cycles in
the sense of \cite{BV2}.

Let $X^p$ be the set of codimension $p$ points in $X$.
For $x\in X^p$ let  $$H^{*}_x(X) \df \limdir{U\subset X}
H^*_{\overline{\{x\}}\cap U}(U).$$ Taking direct limits of (\ref{loc}) over
pairs $Z\subset T$ filtered by codimension and applying the exact couple
method
to the resulting long exact sequence we obtain the coniveau spectral
sequence
$$C^{p,q}_1 =\coprod_{x\in X^p}^{} H^{q+p}_x(X) \implies
H^{p+q}(X)$$
in the abelian category ${\rm  MHS}^{\infty}$ (\cf \cite{BV2}).

Consider $X$ smooth over $\C$. By local purity, we have that
$H^{q+p}_x(X,\Z
(r))\cong H^{q-p}(x,\Z (r-p))$ if $x$ is a codimension $p$ point in $X$,
\ie here
we have set
$$H^*(x) \df\limdir{V {\rm\ open\ } \subset \overline{\{ x\}}} H^*(V).$$
Sheafifying the (limit) sequences (\ref{loc}), we obtain the following
exact sequences of $\Q$-mixed sheaves on $X$:
\begin{equation}\label{shortloc} 0\to
\cH^r_{Z^{p}}\to \coprod_{x\in X^p}^{} x_*
(H^{r-2p}(x)) \to  \cH^{r+1}_{Z^{p+1}}\to 0
\end{equation}
where $\cH^r_{Z^{p}}$ is the $\Q$-mixed sheaf associated to the presheaf
$$U\mapsto  \limdir{{\rm codim}_XZ\geq p} H_{Z\cap U}^r(U).$$
In fact, the claimed short exact sequencs (\ref{shortloc}) are
obtained {\em via}\, the ``locally homologically effaceable'' property (see
\cite[Claim p.  191]{BO}), \ie the following map of sheaves on $X$
 $$\cH^*_{Z^{p+1}} \by{\rm zero}
\cH^*_{Z^p}$$ vanishes for all $p \geq 0$.
\begin{propose} {\em (Arithmetic resolution)} Let $\cH^q_X(\Q(t))$ be the
$\Q$-mixed sheaf defined in (\ref{sheaf}). Assuming $X$ smooth over $\C$
then
\begin{equation}\label{ares}
0\to\cH^q(\Q(t))\to \coprod_{x\in X^0}^{}x_*( H^{q}(x)(t))\to
\coprod_{x\in X^1}^{} x_*( H^{q-1}(x)(t-1))\to \cdots\to
\coprod_{x\in X^q}^{}x_*(\Q (t-q))\to 0
\end{equation}
is a flasque resolution in the category ${\cal MHS}_X$. Therefore, the
coniveau spectral sequence
\begin{equation}\label{coniveau}
C^{p,q}_2 = H^p(X,\cH^q(\Q(t))) \implies
H^{p+q}(X, \Q(t))
\end{equation}
is in the category  ${\rm MHS}^{\infty}$.
\end{propose}
\begin{proof} Follows from construction as sketched above. In fact, all
axioms
stated in \cite[Section~1]{BO} are verified in ${\rm MHS}$ and the results
in \cite[Sections~3-4]{BO} can be obtained in ${\rm MHS}^{\infty}$.
\end{proof}
In particular, consider the presheaf of vector spaces $$U\mapsto F^iH^*(U)
\mbox{\hspace*{1cm} (resp.\ } U\mapsto W_iH^*(U))$$ and the associated
Zariski sheaves $\cF^i\cH^*$ (resp.  $\cW_i\cH^*$) on $X$ filtering the
sheaves $\cH^*(\C)$ (resp.  $\cH^*(\Q)$). These filtrations are defining
the sheaf of mixed Hodge structures $\cH^*_X$ above according to
(\ref{sheaf}).  From Lemma~\ref{abel} (\cf \cite[Theor.1.2.10 and
2.3.5]{D}) the functors $F^i$, $W_i$, $gr_i^W$ and $gr^i_F$ (any $i\in\Z$)
from the category of $\Q$-mixed sheaves to that of ordinary sheaves are
exact.  Applying these functors to the arithmetic resolution (\ref{ares})
we obtain resolutions of $\cF^i\cH^*$ (resp.  $\cW_i\cH^*$) as follows.

\begin{cor}\label{arifilt} The arithmetic resolution (\ref{ares}) yields a
bifiltered
quasi-isomorphism $$(\cH^*,\cF^{\dag},\cW_{\sharp})\by{\simeq}
(\coprod_{x\in X^{\odot}}^{}
x_*H^{*-\odot}(x),\coprod_{x\in X^{\odot}}^{} x_*F^{\dag-\odot},
\coprod_{x\in
X^{\odot}}^{} x_*W_{\sharp-2\odot}),$$ \ie there are flasque resolutions:
$$0\to \cF^i\cH^q\to
\coprod_{x\in X^0}^{} x_*(F^iH^{q}(x))\to \coprod_{x\in X^1}^{} x_*(
F^{i-1}H^{q-1}(x)) \to \cdots$$ and
$$0\to \cW_j\cH^q\to \coprod_{x\in X^0}^{} x_*( W_jH^{q}(x))\to
\coprod_{x\in X^1}^{} x_*(W_{j-2}H^{q-1}(x)) \to \cdots$$
as well as
$$0\to gr^i_{\cF}gr_j^{\cW}\cH^q(\C)\to\coprod_{x\in
X^{0}}^{} x_*(gr^i_{F}gr_j^{W}H^{q}(x))\to\cdots\to
\coprod_{x\in
X^{q}}^{} x_*(gr^{i-q}_{F}gr_{j-2q}^{W}H^{0}(x))\to 0.$$
\end{cor}
Consider the twisted Poincar\'e duality
theory $(F^nH^{*},F^{-m}H_*)$ where the integers $n$ and
$m$ play the role of twisting, \ie we have
$$F^{d-n}H^{2d-k}_Z(X)\cong F^{-n}H_k(Z)$$ for $X$ smooth
of dimension $d$. From the arithmetic resolution of
$\cF^i\cH^q$ in Scholium~\ref{arifilt} we obtain the following:
 \begin{cor}
Assume $X$ smooth and let $i$ be a fixed integer. We then
have a coniveau spectral sequence
\begin{equation}\label{conifilt} F^iC^{p,q}_2 =H^p(X,\cF^i\cH^q)
\implies F^iH^{p+q}(X) \end{equation} where
$H^p(X,\cF^i\cH^q)=0$ if $q<\mbox{min}(i,p)$.
\end{cor}
Concerning the Zariski sheaves $gr^i_{\cF}\cH^q$ and $\bar \cF^i$ we
indeed obtain corresponding coniveau spectral sequences as above.

Note that the spectral sequence $F^iC_2$ can be obtained applying $F^i$ to
the coniveau spectral sequence (\ref{coniveau}).

\begin{remark}{\em \begin{description}
\item[\it i)] Remark that applying $F^i$ (resp.  $W_i$) to the long exact
sequences
(\ref{loc}), taking direct limits over pairs $Z\subset T$ filtered by
codimension and sheafifying, we do obtain the claimed flasque resolutions
of $\cF^i\cH^*$ and $\cW_i\cH^*$ without reference to the category of
$\Q$-mixed sheaves.

\item[\it ii)] Note that for $X$ of dimension $d$, the
fundamental class $\eta_X$ belongs to $W_{-2d}H_{2d}(X)\cap
F^{-d}H_{2d}(X)$ so that ``local purity'' yields the shift by two for the
weight filtration and the shift by one for the Hodge filtration. Therefore
one has to keep care of Tate twists when dealing with arithmetic
resolutions.

\item[\it iii)] Note that for $x\in X^0$ we have $H^{q}(x) = \cH_x^q$ and
there
is a natural projection in ${\cal MHS_X}$ $$\prod_{x\in X}^{}
x_*(\cH_x^q)\to
\coprod_{x\in X^0}^{}x_*(H^{q}(x))$$ which is the identity on $\cH^q$.  The
mixed Hodge structure induced by the arithmetic resolution on $H^*(X,
\cH^q)$ is not the canonical one (which is the one induced by the canonical
$\Q$-mixed flasque resolution) but yields a mixed Hodge structure which is
naturally isomorphic to the canonical one (being induced by a natural
isomorphism in the derived category $\cD^*({\cal MHS}_X)$).
\end{description}}
\end{remark}

\subsection{Coniveau filtration}
Let $X$ be a smooth $\C$-scheme. The coniveau filtration (\cf
\cite{GH}) $N^iH^j(X)$ is a filtration by (mixed) sub-structures of
$H^j(X)$.
This filtration is clearly induced from the coniveau spectral sequence
(\ref{coniveau}) {\em via}\, (\ref{loc}).
Remark that from the coniveau spectral sequence (\ref{coniveau})
$$\gr_N^{i-1}H^j(X) =C^{i-1, j-i+1}_{\infty}$$
which is a substructure of $H^{i-1}(X,\cH^{j-i+1})$ for $i\leq 2$. In fact,
from the arithmetic resolution we have that $C^{p,q}_2=H^p(X,\cH^q(t))=0$
for $p>q$.

\subsubsection*{Case $i=1$}
Let $X$ be a proper smooth $\C$-scheme.
Note that $N^1H^j(X) = \ker (H^j(X)\to H^0(X,\cH^j))=\{\mbox{\ Zariski
locally trivial classes in\ } H^j(X)\}.$ Thus $$\frac{H^j(X, \Q)\cap
F^1H^j(X)}{N^1H^j(X)} = \gr_N^{0}H^j(X)\cap F^1\subseteq H^0(X,\cH^j)\cap
F^1.$$ We remark that $\cH^j/\cF^1$ is the constant sheaf associated to
$H^j(X,\cO_X)$.  Thus $$F^1\cap H^0(X,\cH^j) \cong \ker (H^0(X,\cH^j)\to
H^j(X,\cO_X)).$$ If $j=1$ then $H^1(X)= H^0(X,\cH^1)$ from (\ref{coniveau})
and (\ref{coin}) is trivially an equality.  If $j=2$ then $F^1\cap
H^0(X,\cH^2)=0$ from the exponential sequence.  But for $j=3$ and $X$ the
threefold product of an elliptic curve with itself Grothendieck's argument
in \cite{GH} yields a non-trivial element in $F^1\cap H^0(X,\cH^3)$.

\subsubsection*{Case $i=p$ and $j=2p$}
Let $X$ be a smooth $\C$-scheme.
If $j=2p$ then $N^iH^{2p}(X)(t)=0$ for $i>p$ and
$N^pH^{2p}(X)(t)=C^{p,p}_{\infty}$. Moreover, from (\ref{coniveau})
there is an induced edge map
$$s\ell^{p}_0: H^p(X,\cH^p_X(p))\to H^{2p}(X)(p) $$
which is a map of $\infty$-mixed Hodge structures and whose image is
$N^pH^{2p}(X)(p)$. This equal the image of the classical cycle class map
$c\ell^{p}: CH^p(X)\to H^{2p}(X)(p)$. In fact, by \cite[7.6]{BO},
the cohomology group
$$H^p(X,\cH^p_X(p))\cong {\rm coker} (\coprod_{x\in
X^{p-1}}^{}H^{1}(x) \to \coprod_{x\in X^p}^{} \Z)$$
coincide with $NS^p(X)$, the group of algebraic cycles of
codimension $p$ in $X$ modulo algebraic equivalence. Thus $c\ell^{p}$
factors through
$s\ell^{p}_0$ and the canonical projection (see \cite{BV2}).

Recall that $F^iH^p(X,\cH^q)\df H^p(X,\cF^i\cH^q)\into H^p(X,\cH^q(\C))$
is injective
and $$H^p(X,\cF^p\cH^p)\cong {\rm coker} (\coprod_{x\in
X^{p-1}}^{}F^1H^{1}(x) \to \coprod_{x\in X^p}^{} \C)$$ whence the canonical
map $H^p(X,\cF^p\cH^p)\to NS^p(X)\otimes\C$ is also surjective.  As an
immediate consequence of this fact, \eg from the coniveau spectral sequence
(\ref{conifilt}), we get the following.

\begin{cor}\label{Nero} Let $X$ be a proper smooth $\C$-scheme.
Then $$F^0H^p(X,\cH^p(p)) \df H^p(X,\cF^0\cH^p(p))\cong NS^p(X)\otimes\C$$
and the image of the cycle map is in $H^{2p}(X,\Q (p))\cap
F^0H^{2p}(X, \C (p))$.
\end{cor}
Now $\im c\ell^p_{\sQ}=N^pH^{2p}(X, \Q (p))$ and $H^{2p}(X, \Q (p))\cap
F^0H^{2p}(X, \C (p))$ is equal to $H^{p,p}_{\sQ}$, \ie the sub-structure
of rational $(p,p)$-classes in $H^{2p}(X)$. Note that $H^{p,p}_{\sQ}$
corresponds to the 1-motivic part of $H^{2p}(X)(p)$. For $X$ a smooth
proper $\C$-scheme, the Hodge conjecture then claims that $\im
c\ell^p_{\sQ} = H^{p,p}_{\sQ}$.

In this case $\gr_N^{p-1}H^{2p}(X)(p)=C^{p-1,p+1}_{\infty}$ is a quotient
of $H^{p-1}(X,\cH^{p+1}_X(p))$.  For example: $\gr_N^{1}H^{4}(X)(2)=
H^{1}(X,\cH^{3}_X(2)).$

\begin{question} Is $F^2\cap H^{1}(X,\cH^{3}_X) =0$ ?
\end{question}

\subsubsection*{Case $i=p$ and $j=2p+1$}
Let $X$ be a smooth $\C$-scheme.
If $j=2p+1$ then $N^iH^{2p+1}(X)=0$ for $i>p$ and
$N^pH^{2p+1}(X)(t)=C^{p,p+1}_{\infty}$ which is a quotient of
$H^{p}(X,\cH^{p+1}_X)$, \ie there is an edge map
$$s\ell^{p+1}_{-1}: H^{p}(X,\cH^{p+1}_X(p+1))\to H^{2p+1}(X)(p+1)$$ with
image $N^pH^{2p+1}(X)(p+1)$.  In this case the Grothendieck-Hodge
conjecture
characterize $N^pH^{2p+1}(X)$ as the largest sub-Hodge structure of type
$\{(p,p+1), (p+1,p)\}$.  This is the same as the 1-motivic part of
$H^{2p+1}(X)(p+1)$.

This 1-motivic part yields an abelian variety which is the maximal abelian
subvariety of the intermediate jacobian $J^{p+1}(X)$. On the other hand,
it is easy to see that $N^pH^{2p+1}(X)(p+1)$ yields the algebraic part of
$J^{p+1}(X)$, \ie defined by the images of codimension $p+1$ cycles on $X$
which are algebraically equivalent to zero modulo rational equivalence
(\cf \cite{GH} and \cite{MU2}).

\subsection{Exotic $(1,1)$-classes}
Consider $X$ singular. We briefly explain the Conjecture~\ref{MHC} for
$p=1$.
Moreover we show that there are edge maps generalizing the cycle class maps
constructed in the previous section.

For $X$ a proper irreducible $\C$-scheme, consider the mixed Hodge
structure on $H^{2 + i}(X,\Z )$ modulo torsion.  The extension (\ref{ext})
is the following
\begin{equation}\label{ext1}
0\to H^{1+ i}((H^1)^{\bullet})\to W_2H^{2 + i}(X)/W_0 \to
H^{i}((H^2)^{\bullet})\to 0.
\end{equation}
Since the complex $(H^1)^{\bullet}$ is made of level $1$ mixed Hodge
structures then
$H^{2 + i}(X)^h = H^{2 + i}(X)^e$ in our notation.

If $X$ is nonsingular then $H^{2 + i}(X)$ is pure and there are only two
cases where this extension is non-trivial. In the case $i= -1$ the above
conjecture corresponds to the well known fact that $H_1(\Pic^0(X))= H^1(X,
\Z)$. The case $i = 0$ corresponds to the celebrated theorem by Lefschetz
showing that the subgroup $H^{1,1}_{\sZ}$ of $H^2(X,\Z)$ of cohomology
classes of type $(1,1)$ is generated by $c_1$ of line bundles on $X$.
Since
homological and algebraic equivalences coincide for divisors, the
N\'eron-Severi group $\NS^1 (X)$ coincide with $H^{1,1}_{\sZ}$.
For such a nonsingular variety $X$ we then have $$\NS^1 (X) = F^1\cap
H^2(X,\Z) =
H^{1,1}_{\sZ} = H^1(X,\cH^1_X) = N^1H^2(X,\Z).$$

For $i= -1$ and $X$ possibly singular, the conjecture corresponds to the
fact (proved in \cite{BSAP}) that the abelian variety corresponding to
$\gr_1^WH^1$ is $\ker^0 (\Pic^0(X_0) \to \Pic^0(X_1))$.

For $i = 0$ the Conjecture~\ref{MHC} is quite easily verified by checking
the claimed compatibility of the extension class map. Such a statement then
corresponds to a Lefschetz $(1,1)$-theorem for complete varieties with arbitrary
singularities.

For $i\geq 1$ we may get {\em exotic $(1, 1)$-classes}\, in the higher
cohomology groups $H^{2 + i}(X)$ of an higher dimensional singular variety
$X$. We ignore the geometrical meaning of these exotic $(1, 1)$-classes.
It will be interesting to produce concrete examples. The conjectural
picture is as follows.

Let $\pi : X_{\d}\to X$ be an hypercovering. Let $(H^q(\cH^1))^{\bullet}$
be the complex of $E^{\bullet ,q}_1$-terms of the spectral sequence in
Corollary \ref{zarhodge} for $r=1$. Now $E^{i ,q}_1 =
H^q(X_i,\cH_{X_i}^1)=0$ for $q\geq 2$ (where $X_i$ are the smooth
components of the hypercovering $X_{\d}$ of $X$) and all non-zero terms
are pure Hodge structures: therefore the spectral sequence degenerates at
$E_2$. Thus, from Corollary \ref{zarhodge}, we get an extension
\begin{equation}\label{zarext1}
0\to H^{1+ i}((H^0(\cH^1))^{\bullet})\to \HH^{1 +i}(\Xs, \cH^{1}_{\Xs})
\to H^{i}((H^1(\cH^1))^{\bullet})\to 0
\end{equation}
in the category of mixed $\Q$-Hodge structures.
We have $H^{i+1}((H^0(\cH^1))^{\bullet})=
H^{i+1}((H^1)^{\bullet})= \gr_1^WH^{2 + i}$
and $H^{i}((H^1(\cH^1))^{\bullet})= H^{i}((NS)^{\bullet}) =
H^{i}((N^1H^{2})^{\bullet})$.

Moreover, from the local-to-global spectral sequence in Claim~\ref{l2g}
and cohomological descent we get the following edge map
$$s\ell^{1+i}:\HH^{1 +i}(\Xs, \cH^{1}_{\Xs})\to W_2H^{2 + i}(X)/W_0.$$
In fact, first observe that $W_0H^{2 + i}(X) = \HH^{2 +i}(\Xs,
\cH^{0}_{\Xs}).$ From (\ref{zarext1}) above we then see that $W_0\HH^{1
+i}(\Xs, \cH^{1}_{\Xs})=0.$ The map $s\ell^{1+i}$ is then easily obtained
as an edge homomorphism of the cited local-to-global spectral sequence and
weight arguments. This cycle map will fit in a diagram
$$\begin{array}{ccccccc}
0\to & H^{1+ i}((H^1)^{\bullet})&\to & W_2H^{2 + i}(X)/W_0 &\to
&H^{i}((H^2)^{\bullet})&\to 0\\
&\uparrow\veq& &\uparrow {\scriptsize s\ell^{1+i}}& &\uparrow &\\
0\to& H^{1+ i}((H^0(\cH^1))^{\bullet})&\to &\HH^{1 +i}(\Xs,
\cH^{1}_{\Xs})& \to & H^{i}((H^1(\cH^1))^{\bullet})&\to 0
\end{array}$$
mapping the extension (\ref{zarext1}) to (\ref{ext1}).
\begin{cor} The image of the map $s\ell^{1+i}$ is $H^{2 + i}(X)^e$.
\end{cor}

Following \cite{BSAP} consider the simplicial sheaf $\cO_{\Xs}^*$ and the
corresponding Zariski cohomology groups
$\HH^{1 +i}(\Xs,\cO_{\Xs}^*)$.
Since the components of $\Xs$ are smooth, the canonical spectral sequence
$$E^{p ,q}_1 = H^q(X_p, \cO^*_{X_p}) \implies \HH^{p+ q}(\Xs,\cO_{\Xs}^*)$$
yields a long exact sequence
$$\cdots \to H^{1+ i}((H^0(\cO^*))^{\bullet})\to \HH^{1
+i}(\Xs,\cO_{\Xs}^*) \to H^{i}((\Pic))^{\bullet})\longby{d^{i}} H^{2+
i}((H^0(\cO^*))^{\bullet}) \to \cdots $$
According to \cite{BSAP} (see the construction in \cite{BRS}) we may regard
$\HH^{1+i}(\Xs,\cO_{\Xs}^*)$ as the group of $k$-points of a group scheme
whose connected component of the identity yields a semi-abelian variety
$$0 \to H^{1+ i}((H^0(\cO^*))^{\bullet})/\sigma\to \HH^{1
+i}(\Xs,\cO_{\Xs}^*)^0 \to H^{i}((\Pic^0))^{\bullet})^0\to 0$$
where $\sigma$ is a finite group. The Hodge realization of the so obtained
isogeny 1-motive is
$$T_{\rm Hodge} ([0\to \HH^{1+i}(\Xs,\cO_{\Xs}^*)^0]_{\sQ}) =
W_1H^{1+i}(X,\Q)(1).$$
This last claim is clearly related to Deligne's conjecture
\cite[10.4.1]{D}. For $i=-1, 0$ this is actually proven in \cite{BSAP} and
for all $i$ in \cite{BRS}.

Recall the existence of a canonical map of sheaves $c_1: \cO_{\Xs}^*\to
\cH^{1}_{\Xs}$ yielding a map $$c_1: \HH^{1 +i}(\Xs,\cO_{\Xs}^*)\to \HH^{1
+i}(\Xs, \cH^{1}_{\Xs}).$$ By composing $s\ell^{1+i}$ and $c_1$ we then
obtain a cycle map $$\HH^{1 +i}(\Xs,\cO_{\Xs}^*)\to W_2H^{2 + i}(X)/W_0.$$
We may regard the image of this cycle map as the discrete part of $\HH^{1
+i}(\Xs,\cO_{\Xs}^*)$.  Over $\Q$, it is clearly equal to $F^1\cap H^{2 +
i}(X,\Q)$.  The reader can easily check that this is the case, \eg $\HH^{2
+ i}(\Xs,\cO_{\Xs}^*[-1])$ coincides with Deligne-Beilinson
cohomology (see \cite[5.4]{BE}).  In general, we may expect the following
picture for cycle maps.

\subsection{$\cK$-cohomology and motivic cohomology}
Let $\Xs$ be a smooth simplicial scheme.
Consider the local-to-global spectral sequence in Claim~\ref{l2g}. For a
fixed $p$ we then obtain a spectral sequence
$$W_{2p}\HH^q(X_{\d}, \cH_{\Xs}^r)/W_{2p-2}\implies
W_{2p}\HH^{q+r}(X_{\d},\Q)/W_{2p-2}.$$
The sheaf $\cH_{\Xs}^r$ has weights $\leq 2r$ and so the mixed Hodge
structure
on its cohomology has weights $\leq 2r$.
Thus $W_{2p}\HH^q(X_{\d}, \cH_{\Xs}^p) = \HH^q(X_{\d}, \cH_{\Xs}^p)$ and
$W_{2p}\HH^q(X_{\d}, \cH_{\Xs}^r)/W_{2p-2}=0$ if $r<p$. Thus, there is an
edge map
\begin{equation}\label{sedge}
s\ell^{p+i}:\HH^{p+i}(X_{\d}, \cH_{\Xs}^p)/W_{2p-2}\to
W_{2p}\HH^{2p+i}(\Xs)/W_{2p-2}.
\end{equation}
Note that if $\Xs = X$ is constant then $s\ell^{p+0}=s\ell^{p}_0$ and
$s\ell^{p-1}=s\ell^{p}_{-1}$ in the notation of Section~4.2, modulo
$W_{2p-2}$.
\begin{cor} The image of the edge map $s\ell^{p+i}$ is $H^{2p
+ i}(X)^h$.\end{cor}

Consider Quillen's higher $K-$theory. Consider Zariski sheaves associated
to Quillen's $K$-functors. The $\cK$-cohomology groups are $\HH^*(X_{\d},
\cK_p)$ (as usual we consider Zariski simplicial sheaves $\cK_p$). Local
higher Chern classes give us maps of simplicial sheaves $c_p : \cK_p \to
\cH_{\Xs}^p(p)$ for each $p\geq 0$ (\cf \cite{BV2}). We thus obtain a map
\begin{equation}
c_p : \HH^{p+i}(X_{\d}, \cK_p) \to \HH^{p+i}(X_{\d}, \cH_{\Xs}^p(p)).
\end{equation}
Note that in the canonical spectral sequence $$E^{s ,t}_1 = H^t(X_s,
\cK_{p}) \implies \HH^{s+t}(\Xs,\cK_{p})$$ we have $ H^t(X_s, \cK_{p}) =
0$ if $t >p$.  The same hold for the sheaf $\cH_{\Xs}^p$. We then have
a commutative square
$$\begin{array}{ccc}
\HH^{p+i}(X_{\d}, \cK_p)& \by{c_p} &\HH^{p+i}(X_{\d}, \cH_{\Xs}^p(p))\\
\downarrow & &\downarrow\\
H^{i}((CH^p)^{\bullet})&\to & H^{i}((NS^p)^{\bullet}).
\end{array}$$
Thus the image of $c_p$ in $H^{i}((NS^p)^{\bullet})$ is clearly contained
in the kernel of the map $\lambda_a^i$ defined in (\ref{bound}). We also
have the following commutative square
$$\begin{array}{ccc}
\HH^{p+i}(X_{\d}, \cH_{\Xs}^p)/W_{2p-2}& \by{s\ell^{p+i}} &
W_{2p}\HH^{2p+i}(\Xs)/W_{2p-2}\\
\downarrow & &\downarrow\\
H^{i}((NS^p)^{\bullet})&\to & \gr_{2p}\HH^{2p+i}(\Xs).
\end{array}$$
Composing $s\ell^{p+i}$ and $c_p$ above we then obtain a simplicial cycle
map \begin{equation}\label{scycle}
c\ell^{p+i}: \HH^{p+i}(X_{\d}, \cK_p) \to W_{2p}H^{2p+i}(\Xs)/W_{2p-2}.
\end{equation}
Let $X$ be a proper $\C$-scheme and let $\Xs\to X$ be a universal cohomological
descent morphism. By descent,
$\HH^{*}(\Xs)\cong H^{*}(X)$ as mixed Hodge structures. Let $H$ denote the
mixed Hodge structure on $H^{2p+i}(X, \Z)/({\rm torsion})$. Let $F^p$ denote
the Hodge filtration. Note that
$$ F^p\cap H_{\sZ} =
\Hom_{\rm MHS} (\Z(-p), H) =
\Hom_{\rm MHS} (\Z(-p), W_{2p}H).$$
Moreover
$$\Hom_{\rm MHS} (\Z(-p), W_{2p}H) \subseteq
\Hom_{\rm MHS} (\Z(-p), W_{2p}H/W_{2p-2}H) = F^p \cap
H^e_{\sZ}.$$ Thus
$$F^p\cap H_{\sZ}\subseteq  F^p \cap H^e_{\sZ}= F^p \cap H^h_{\sZ} =
\ker (H^{p,p}_{\sZ}\by{e^p} J^p(H)).$$
Therefore we have the following natural question.
\begin{question} \label{cycleimage} Let $X$ be a proper $\C$-scheme and
let $\pi :\Xs\to X$ be a proper smooth hypercovering.  Is $F^p \cap
H^{2p+i}(X, \Q)$ the image of the cycle class map $c\ell^{p+i}$ in
(\ref{scycle}) ?
\end{question}

Bloch's counterexample answer this question in the negative (see Section~5.1 below).
Note that here we actually deal with the simplicial scheme $\Xs$ (not just $X$)
as, \eg in \cite{BS} it is shown that $F^1 \cap H^{2}(X, \Z)$ can be larger
than the image of $\Pic (X)$, if $X$ is singular.  However, $F^1 \cap
H^{2}(X, \Z)$ is the image of the $\Pic$ of any hypercovering of $X$.
However, $F^2 \cap H^{4}(X, \Q)$ is larger than the image of $\HH^2(\Xs,\cK_2)$ if $X$ is
the singular 3-fold in Section~5.1 below.

Let's then consider the case $p =2$ in the above.  In this case we have that
$H^q(X_i,\cH^2_{X_i})$ is purely of weight $q+2$ by the coniveau spectral
sequence (\ref{conifilt}).  Thus the canonical spectral sequence in
Corollary~\ref{zarhodge} degenerates yielding the following extension
\begin{equation}\label{zarext2}
0\to H^{1+ i}((H^1(\cH^2))^{\bullet})\to \HH^{2+i}(\Xs, \cH^{2}_{\Xs})/W_2
\to H^{i}((NS^2)^{\bullet})\to 0.
\end{equation}
Note that $H^1(X_i,\cH^2_{X_i})= N^1H^3(X_i)$ and $W_2\HH^{j}(\Xs,
\cH^{2}_{\Xs})= H^{j}((H^0(\cH^2))^{\bullet})$. The map in (\ref{sedge})
is mapping the extension (\ref{zarext2}) to the following canonical
extension
$$0\to H^{1+ i}((H^3)^{\bullet})\to  W_4H^{4 + i}(X)/W_2 \to
H^{i}((H^4)^{\bullet})\to 0.$$
According to Conjectures~\ref{GHC}--\ref{MHC} we may expect that the image of
$\HH^{2+i}(\Xs, \cH^{2}_{\Xs})/W_2$ under this map is $H^{4 + i}(X)^h$.

Finally, making use of the triangulated category of motives (see
\cite{V} and \cite{LM}) let $H^*_m(X,\Q (\cdot))$ denote the motivic cohomology
of the proper $\C$-scheme $X$.  Since motivic cohomology is universal we may get
a canonical map $H^{2p+i}_m(X,\Q (p))\to H^{2p+i}(X,\Q (p))$ compatibly with the
weight filtrations. This map will factors trhough  Beilinson's absolute Hodge
cohomology \cite{BE}. However, in general, its image will not be larger than
$c\ell^{p+i}$ in (\ref{scycle}), \ie smaller than the rational part of
$F^pH^{2p+i}(X,\C)$.  In fact, we can see that only Beilinson's absolute Hodge
cohomology (or Deligne-Beilinson cohomology) would have image equal to $F^2\cap
H^{4}(X,\Q)$ if $X$ is the singular 3-fold in Bloch's counterexample below.

\section{Examples}
We finally discuss a couple of examples where one can test the conjectures.

\subsection{Bloch's example}
We now consider Bloch's example explained in a letter to U. Jannsen,
reproduced in the Appendix~A of \cite{JA} (see also
Appendix~A.I in \cite{LW}). This example, originally requested by Mumford,
is a counterexample to a naive extension of the cohomological Hodge
conjecture to the singular case. Moreover (as indicated by Bloch's Remark~1 in
\cite[Appendix~A]{JA}) it shows that no cohomological invariants of algebraic
varieties, that agree with Chow groups of non-singular varieties, can provide
all Hodge cycles for singular varieties.

Let $P$ be the blow-up of $\P^3$ at a point $x$ in $S_0\subset \P^3$ a smooth
hypersurface of degree $\geq 4$ over $\bar{\Q}$. The point $x$ is assumed
$\bar{\Q}$-generic. Let $S$ be the blow-up of $S_0$ at $x$ over $\C$. Thus
$S\subset P$ and $H^3(S,\Q)=0$.

Let $X$ be the gluing of two copies of $P$ along $S$, \ie the
singular projective variety defined as the pushout
\[\begin{array}{ccc}
S\coprod S&\by{i\coprod i}& P\coprod P\\
c\downarrow\quad&&\quad\downarrow f\\
S&\by{j}&X
\end{array}\]
Such a Mayer-Vietoris diagram always defines a cohomological descent
morphism $\Xs \to X$ (\eg a distinguished (semi)simplicial resolution in
the sense of Carlson \cite[\S 3 and \S 13]{CA}).

Thus we obtain a short exact sequence $$0 \to H^4(X,\Q
(2)) \to H^4(P,\Q(2))^{\oplus 2} \to H^4(S,\Q(2)) \to 0$$
where $CH^2(P)_{\sQ} \cong H^4(P,\Q(2))$.
Thus $H^4(X, \Z)$ has rank 3, is purely of type $(2,2)$
and the Hodge 1-motive is $[H^4(X,\Q (2))  \to 0 ]$. From (\ref{zarext2})
we obtain $\HH^{2}(\Xs, \cH^{2}_{\Xs})/W_2 = H^0((\NS^2)^{\bullet})) =H^4(X,\Q
(2)).$

Since the Albanese of $S$ vanishes, the algebraically defined Hodge 1-motive
is given by $H^0((\NS^2)^{\bullet})) =\ker (NS^2(P)^{\oplus 2} \to NS^2(S))$ and
we clearly have that
$$[H^0((\NS^2)^{\bullet})) \to 0]\cong [H^4(X,\Q (2))  \to 0 ]$$
as predicted by Conjecture~\ref{MHC}. However
$H^0((CH^2)^{\bullet})) = \ker (CH^2(P)^{\oplus 2} \to CH^2(S))$ has rank 2, and
it is strictly smaller than $H^4(X,\Q (2))$, as Bloch's observed. Moreover, from
the above we may regard  $\HH^2(\Xs,\cK_2)_{\sQ}$ mapping to both
$H^0((CH^2)^{\bullet}))$ and $H^4(X,\Q (2))$. Then $c\ell^{2}$ in
(\ref{scycle}) is not surjective because $H^0((CH^2)^{\bullet}))\neq H^4(X,\Q
(2))$. The same argument applies to motivic cohomology $H^4_m(X,\Q (2))$. In
fact, if $Y$ is a smooth variety $H^4_m(Y,\Q (2))\cong CH^2(Y)_{\sQ}$. However,
Beilinson's absolute Hodge cohomology is
$H^4_{\cD}(X,\Z (2))\otimes \Q\cong H^4(X,\Q (2))$. Thus, Beilinson's formulation
of the Hodge conjecture in \cite[\S 6]{BE} doesn't hold in the singular case.

Note that in this example, all Hodge classes are involved, as the Hodge structure
is  pure.

\subsection{Srinivas example}
The following example has been produced by Srinivas upon author's request.
It is similar to Bloch's example however, in this example, the space of
``Hodge cycles" is strictly smaller than $H^{2,2}_{\sQ}$.

Let $Y$ be a smooth projective complex 4-fold with
$H^1(Y,\Z)=H^3(Y,\Z)=0$, and with an algebraic cycle $\alpha\in CH^2(Y)$
whose singular cohomology class $\bar{\alpha}\in H^4(Y,\Q)$ is a non-zero
primitive class. For example, $Y$ could be a smooth quadric hypersurface
in $\P^5$, and $\alpha\in CH^2(Y)$ the difference of the classes of two
planes, taken from the two distinct connected families of planes in $Y$.
Let $Z$ be a general hypersurface section of $Y$ of any fixed degree
$d$ such that $H^{3,0}(Z)\neq 0$ (this holds for any large enough degree
$d$; for example, if $Y$ is a quadric then we may take $Z=Y\cap H$ to be
the intersection with a general hypersurface $H$ of any degree $\geq 3$).

Then $Z$ is a smooth projective 3-fold, and by the theorem of Griffiths,
if $i:Z\to Y$ is the inclusion, then $i^*\alpha\in CH^2(Z)$ is
homologically trivial, but no non-zero multiple of $i^*\alpha$
is algebraically equivalent to 0. In fact, $H^3(Z,\Q)$ has no proper
Hodge substructures, and so (because $H^{3,0}(Z)\neq 0$) the Abel-Jacobi
map vanishes on the group $CH^2(Z)_{\rm alg}$ of cycle classes
algebraically equivalent to 0; on the other hand, the Abel-Jacobi image
of $i^*\alpha$ is non-torsion.

Now let $X$ be the singular projective variety defined as a push-out
\[\begin{array}{ccc}
Z\coprod Z&\by{i\coprod i}& Y\coprod Y\\
c\downarrow\quad&&\quad\downarrow f\\
Z&\by{j}&X
\end{array}\]
so that $X$ is obtained by gluing two copies of $Y$ along $Z$.

Consider the simplicial scheme $X_{\d}$ obtained as above (\eg the \v{C}ech
hypercovering of $X$, with $X_0\to X$ taken to be the quotient map
$f:Y\coprod Y\to X$). Then $H^*(X_{\d},\Q)\cong H^*(X,\Q)$ as mixed Hodge
structures, and we have an exact sequence of mixed Hodge structures (of which all
terms except $H^4(X,\Z)$ are in fact pure)
\[0\to H^3(Z,\Z)\to H^4(X,\Z)\to H^4(Y,\Z)^{\oplus 2}\by{s} H^4(Z,\Z)\]
where $s(a,b)=i^*a-i^*b$. Then $(\bar{\alpha},0)$ and $(0,\bar{\alpha})$
are linearly independent elements of $\ker s$, since $i^*\bar{\alpha}=0$
in $H^4(Z,\Q)$ (this is essentially the definition of $\bar{\alpha}$
being a primitive cohomology class).

In this situation, the group of Hodge classes in $H^4(X,\Q)/W_3$ is
non-trivial, but since $H^3(Z,\Q)$ has no non-trivial sub-Hodge
structures, the intermediate Jacobian $J^2(Z)$ has no non-trivial abelian
subvariety. The extension of Hodge structures determined by
the Hodge classes is not split; for example the extension class of the
pullback of
\[0\to H^3(Z,\Z)\to H^4(X,\Z)\to \ker s\to 0\]
under $\Z(-2)\to \ker s$ determined by $(\bar{\alpha},0)$ is (up to sign)
the Abel-Jacobi image of $i^*\alpha$, which is non-torsion.
Here $\rank (\ker s)=3$, so we get an extension class map
$\Z^3\to J^2(Z)$; one checks that the image has rank 1, generated by the image
of $(\bar{\alpha},0)$ (or equivalently by the image of $(0,\bar{\alpha})$).

So the lattice for the corresponding Hodge 1-motive is, by
definition $$\ker (\Z^3\to J^2(Z))=F^2\cap H^4(X,\Z),$$ which is strictly smaller
than the lattice of all Hodge classes in $H^4(X,\Z)$. Moreover, since
$H^1(Z,\cH^2)= N^1H^3(Z)=0$, from the extension  (\ref{zarext2}) we obtain
$\HH^{2}(\Xs, \cH^{2}_{\Xs})/W_2 \cong H^{0}((NS^2)^{\bullet})$.

Finally, the cycle $\alpha\in CH^2(Y)$ projects to a cycle in $\NS^2(Y)$
which restricts to a non-zero class $i^*\alpha\in NS^2(Z)_{\sQ}$ by construction.
Since $CH^2(Z)_{\rm ab}=CH^2(Z)_{\rm alg}$ then the algebraically defined
1-motive is given by the image of $H^{0}((NS^2)^{\bullet})=\ker (NS^2(Y)^{\oplus
2} \to NS^2(Z))$ in $H^4(X,\Z)$, providing generators for $\ker (\Z^3\to J^2(Z))$
as claimed in Conjecture~\ref{MHC}.

\vspace*{2cm}
{\sc E-Mail:} {\tt barbieri@dmmm.uniroma1.it}\\
\begin{flushright}
\noindent{\em Dipartimento di Metodi e Modelli Matematici\\ Universit\`a
degli
Studi di Roma ``La Sapienza''\\ Via A. Scarpa, 16\\I-00161 --- Roma
(Italy)}
\end{flushright}
\end{document}